\documentclass[a4paper,10pt]{amsart}

\usepackage[english]{babel}
\usepackage[latin1]{inputenc}
\usepackage{amsmath}
\usepackage{amsthm}
\usepackage{amsfonts}
\usepackage{amssymb, latexsym, graphics, showidx}

\numberwithin{equation}{section}
\newtheorem{de}{Definition}[section]
\newtheorem{thm}{Theorem}[section]
\newtheorem{rem}[thm]{Remark}
\newtheorem{cor}[thm]{Corollary}
\newtheorem{prop}[thm]{Proposition}
\newtheorem{lem}[thm]{Lemma}

\renewcommand{\dim}{\noindent\textbf{Proof.} }
\newcommand{\dims}{\noindent\textbf{Proof} }
\newcommand{\finedim}{{\unskip\nobreak\hfil\penalty50
   \hskip2em\hbox{}\nobreak\hfil\mbox{$\Box$ \qquad}
   \parfillskip=0pt \finalhyphendemerits=0\par\medskip}}
\newcommand{\R}{\mathbb{R}}
\newcommand{\N}{\mathbb{N}}
\newcommand{\Om}{\Omega}
\newcommand{\Oms}{\overline{\Omega}}
\newcommand{\lam}{\lambda}
\newcommand{\al}{\alpha}
\newcommand{\lams}{\overline{\lambda}}
\newcommand{\xs}{\overline{x}}
\newcommand{\ys}{\overline{y}}
\newcommand{\zs}{\overline{z}}
\newcommand{\tr}{\text{tr}}
\newcommand{\di}{\Delta_\infty}
\newcommand{\ep}{\epsilon}

\title[The Neumann eigenvalue problem for the $\infty$-Laplacian]
{The principal eigenvalue of the $\infty$-Laplacian with the
Neumann boundary condition}
\author{Stefania Patrizi}
\address{SAPIENZA Universit\`a di Roma, Dipartimento di Matematica,
Piazzale A.~Moro 2, I-00185 Roma, Italy}
\email{patrizi@mat.uniroma1.it}
\begin{document}
\begin{abstract}
We prove the existence of a principal eigenvalue associated to the
$\infty$-Laplacian plus lower order terms and the Neumann boundary
condition in a bounded smooth domain. As an application  we get
uniqueness and existence results for the Neumann problem and a
decay estimate for viscosity solutions of the Neumann evolution
problem.
\end{abstract}

\maketitle

\section{Introduction}
In this paper we study the maximum principle, the principal
eigenvalue, regularity, existence and uniqueness for viscosity
solutions of the Neumann boundary value problem

\begin{equation}\label{sisintro}
\begin{cases}
\di u +b(x)\cdot Du + (c(x)+\lam)u= g(x)  & \text{in} \quad\Om \\
 \frac{\partial u}{\partial\overrightarrow{n}} = 0 & \text{on} \quad\partial\Om, \\
 \end{cases}
 \end{equation}
 where $\Om$ is a bounded smooth domain,
$\overrightarrow{n}(x)$ is the exterior normal to the domain $\Om$
at $x$, $b$, $c$ and $g$ are continuous functions on $\Oms$,
$\lam\in\R$ and
\begin{equation}\label{inflapl}\di u=\langle
D^2u\frac{Du}{|Du|},\frac{Du}{|Du|}\rangle,\end{equation} for
$u\in C^2(\Om)$, is the 1-homogeneous version of the
$\infty$-Laplacian.

The  $\infty$-Laplacian, which arises from the optimal Lipschitz
extension problem, see \cite{acj}, appears also in the
Monge-Kantorovich mass transfer problem, see \cite{eg}, and
recently, some authors have introduced a game theoretic
interpretation of it, see \cite{pssw}.

We define and investigate the properties of the principal
eigenvalue of the ope\-ra\-tor $$-(\di +b(x)\cdot D+c(x)),$$ with
the Neumann boundary condition and as an application, we get
existence and uniqueness results for \eqref{sisintro} and a decay
estimate for the solution of the associated evolution problem.

In their famous work \cite{bnv}, Berestycki, Nirenberg and
Varadhan defined the principal eigenvalue $\lam_1$ of a general
linear uniformly elliptic operator $-L$ where
$$L[u]=\text{tr}(A(x)D^2u)+b(x)\cdot Du+c(x)u,$$
in a bounded domain $\Om$, as the supremum of those $\lam$ for
which there exists a positive supersolution of $L[u]+\lam u=0$. In
that paper, they showed that $\lam_1$ is the first eigenvalue of
$L$, i.e., for any eigenvalue $\lam\neq \lam_1$,
Re$(\lam)>\lam_1$; moreover $\lam_1$ can be characterized as the
supremum of those $\lam$ for which the operator $L+\lam I$
satisfies the maximum principle, i.e., for any $\lam<\lam_1$, if
$u$ is a subsolution of $L[u]+\lam u=0$ and $u\leq 0$ on
$\partial\Om$ then $u\leq 0$ in $\Om$. They established other
properties of the first eigenvalue, such as simplicity and
stability.

In view of its relation with the maximum and the comparison
principles, the concept of principal eigenvalue has been extended
to nonlinear operators to study the associated boundary value
problems. That has been done for the variational operators, such
that the p-Laplacian, through the method of minimization of the so
called nonlinear Rayleigh quotient,  see e.g. \cite{an} and
\cite{li2}. That method uses heavily the variational structure and
cannot be applied to operators which have not this property. An
important step in the study of the eigenvalue problem for general
nonlinear operators was made by Lions in \cite{l}. In that paper,
using probabilistic and analytical methods, he showed the
existence of principal eigenvalues for the uniformly elliptic
Hamilton-Jacobi-Bellman operator.  Very recently, many authors,
inspired by  \cite{bnv}, have de\-ve\-lo\-ped an eigenvalue theory
for fully nonlinear operators which are non-variational. The
Pucci's extremal operators have been treated by Quaas \cite{q} and
Busca, Esteban and Quaas \cite{beq}. Their results have been
extended to more general fully nonlinear convex uniformly elliptic
operators in \cite{qs} by Quaas and Sirakov. See also the work of
Ishii and Yoshimura \cite{iy} for non-convex operators.

Issues similar to those of this paper have been studied  by
Birindelli and Demengel in \cite{bd} and the author of this note
in \cite{p1} where respectively the Dirichlet and the Neumann
eigenvalue problem is treated for degenerate or singular elliptic
operators $F(x,Du,D^2u)$ plus lower order terms. In these papers,
among other assumptions, $F$  is required to satisfied
\begin{equation}\label{unifellitt} a|p|^\al \text{tr}N\leq F(x,p,M+N)-F(x,p,M)\leq A|p|^\al \text{tr}N,\end{equation}
with $\al>-1$, for $x\in\overline{\Om},\, p\in
\R^N\setminus\{0\},$ and $M, N$ symmetric matrices with $N\geq 0$.
Typical examples are given by $|Du|^\al\mathcal{M}_{a,A}(D^2u)$,
$\al>-1$, where $\mathcal{M}_{a,A}(D^2u)$ is one of the Pucci's
operator, the p-Laplacian and some non-variational generalizations
of it. Because of its strong degeneracy, the $\infty$-Laplacian
does not satisfy \eqref{unifellitt}, so it is not covered by
\cite{bd} or \cite{p1}.

 The existence of a principal eigenvalue
defined as in \cite{bnv} for the $\infty$-Laplacian with the
Dirichlet boundary condition has been treated by Juutinen in
\cite{j} together with many other questions. We want to mention
that there exists also a different approach to investigate the
eigenvalue problem for \eqref{inflapl} which consists in studying
the asymptotic behavior, as $p\rightarrow \infty$, of the
p-Laplacian eigenvalue equation, see \cite{jlm} and \cite{gmpr}.
This second method uses the variational formulation of the
approximate problems and leads to a different limit eigenvalue
problem, see \cite{j}.

Following the ideas of \cite{bnv}, we define the principal
eigenvalue as
\begin{equation}\begin{split}\label{lams}\lams:=\sup\{&\lam\in
\R\;|\;\exists\,v>0 \text{ on $\Oms$ bounded viscosity
supersolution of }\\& \di v +b(x)\cdot Dv + (c(x)+\lam)v= 0 \text{
in } \Om,\, \frac{\partial v}{\partial\overrightarrow{n}} =
0\text{ on }\partial \Om \}.
\end{split}\end{equation}
$\lams$ is well defined since the above set is not empty; indeed,
$-|c|_\infty$ belongs to it, being $v(x)\equiv1$ a corresponding
supersolution. Furthermore it is an interval because if $\lam$
belongs to it then so does any $\lam'<\lam$.

One of the scope of this work is to prove that $\lams$ is an
"eigenvalue" for $-(\di +b(x)\cdot D+c(x))$ which admits a
positive "eigenfunction". As in the linear case it can be
characterized as the supremum of those $\lam$ for which $\di
+b(x)\cdot D+c(x)+\lam$ with the Neumann boundary condition
satisfies the maximum principle. As a consequence, $\lams$ is the
least "eigenvalue", i.e., the least number for which there exists
a non-zero solution of
\begin{equation*}
\begin{cases}
\di u +b(x)\cdot Du + (c(x)+\lam)u= 0  & \text{in} \quad\Om \\
 \frac{\partial u}{\partial\overrightarrow{n}} = 0 & \text{on} \quad\partial\Om. \\
 \end{cases}
 \end{equation*}
These results are applied to obtain existence and uniqueness for
the boundary value problem \eqref{sisintro}.

Remark that since $\di (-u)=-\di u$, $\lams$ can be defined also
in the following way
\begin{equation}\begin{split}\label{lamsotto}\lams=\sup\{&\lam\in
\R\;|\;\exists\,u<0 \text{ on $\Oms$ bounded viscosity subsolution
of }\\& \di u +b(x)\cdot Du + (c(x)+\lam)u= 0 \text{ in } \Om,\,
\frac{\partial u}{\partial\overrightarrow{n}} = 0\text{ on
}\partial \Om \}.
\end{split}\end{equation} For a fully nonlinear operator, $\lams$ defined as in \eqref{lams} may be different
from the quantity defined as in \eqref{lamsotto}, see \cite{p2}.

The paper is organized as follows. In the next section we give
assumptions and precise the concept of solution adopted. In
Section 3 we establish a Lipschitz re\-gu\-la\-rity result for
viscosity solutions of \eqref{sisintro}. Section 4 is devoted to
the maximum principle for subsolutions of \eqref{sisintro}.  In
Section 4.1 we show that it holds (even for more general boundary
conditions) for $\di +b(x)\cdot D+c(x)$ if $c(x)\leq 0$ and
$c\not\equiv 0$, see Theorem \ref{macpcless0}. One of the main
result of the paper is that the maximum principle holds for $\di
+b(x)\cdot D+c(x)+\lam $ for any $\lam<\lams$, as we show in
Theorem \ref{maxpneum} of Section 4.2. In particular it holds for
$\di +b(x)\cdot D+c(x)$ if $\lams>0$. Following the example given
in \cite{p1} we show that the result of Theorem \ref{maxpneum} is
stronger than that of Theorem \ref{macpcless0}, i.e., that there
exist some functions $c(x)$ changing sign in $\Om$ for which the
principal eigenvalue of $\di +b(x)\cdot D+c(x)$ is positive and
then for which the maximum principle holds.

In Section 5 we show some existence and comparison theorems. In
particular, we prove that the Neumann pro\-blem \eqref{sisintro}
is solvable for any right-hand side if $\lam<\lams$.

Finally, in Section 6 we prove a decay estimate for solutions of
the Neumann evolution problem.
\section{Assumptions and definitions}
We denote by $\emph{S(N)}$ the space of symmetric matrices on
$\R^N$ equipped with the usual ordering and we fix the norm
$\|X\|$ in $\emph{S(N)}$ by setting
$\|X\|=\sup\{|X\xi|\,|\,\xi\in\R^N,\,|\xi|\leq1\}=\sup\{|\lam|\,:\,
\lam \text{ is an eigenvalue of }X\}.$

Let $\sigma:\R^N\rightarrow \emph{S(N)}$ be the function defined
by
$$\sigma(p):=\frac{p\otimes p}{|p|^2}.$$
The $\infty$-Laplacian can be written as $$\di u=\tr
(\sigma(Du)D^2u),$$ for any $u\in C^2(\Om)$.

It easy to check that $\sigma$ has the following properties:
\begin{itemize}
    \item $\sigma(p)$ is homogeneous of order 0,
i.e., for any $\al\in \R$ and $p\in\R^N$
$$\sigma(\al p)=\sigma(p);$$
    \item For all $p\in\R^N$ $$0\leq
\sigma(p)\leq I,$$ where $I$ is the identity matrix in $\R^N$;
    \item $\sigma(p)$ is idempotent, i.e.,
    $$(\sigma(p))^2=\sigma(p);$$
    \item For any
$p\in\R^N\setminus\{0\}$ and $p_0\in\R^n$ with $|p_0|\leq
\frac{|p|}{2}$
\begin{equation}\label{infp-infq}
\tr\left[(\sigma(p+p_0)-\sigma(p))^2\right]\leq 8
\frac{|p_0|^2}{|p|^2}.
\end{equation}
\end{itemize}

The domain $\Om$ is supposed to be bounded and of class $C^2$. In
particular, it satisfies the interior sphere condition and the
uniform exterior sphere condition, i.e.
\begin{enumerate}
\item [($\Om1$)]For
 each $x\in\partial \Om$ there exist $R>0$ and $y\in \Om$ for
 which $|x-y|=R$ and $B(y,R)\subset \Om$.
\item[($\Om2$)] There exists $r>0$ such that  $B(x+r\overrightarrow{n}(x),r)\cap \Om =\emptyset$  for any $x\in \partial\Om.$\end{enumerate}
From the property ($\Om2$) it follows that
\begin{equation}\label{sferaest}
\langle
y-x,\overrightarrow{n}(x)\rangle\leq\frac{1}{2r}|y-x|^2\quad\text{for
}x\in\partial\Om \text{ and } y\in\Oms.\end{equation} Moreover,
the $C^2$-regularity of $\Om$ implies the existence of a
neighborhood of $\partial\Om$ in $\Oms$ on which the distance from
the boundary
 $$d(x):=\inf\{|x-y|, y\in\partial\Om\},\quad x\in\Oms$$ is of class $C^2$. We still denote by $d$ a
$C^2$ extension of the distance function to the whole $\Oms$.
Without loss of generality we can assume that $|Dd(x)|\leq1$ on
$\Oms$.

We adopt the notion of viscosity solution for \eqref{sisintro}
given in \cite{bd} for singular elliptic operators, in which is
required to test only with test functions which have gradient
different from zero.

 We denote by $USC(\overline{\Om})$
the set of upper semicontinuous functions on $\overline{\Om}$ and
by $LSC(\overline{\Om})$ the set of lower semicontinuous functions
on $\overline{\Om}$. Let $g:\Oms\rightarrow \R$ and
$B:\partial\Om\times\R\times \R^N\rightarrow\R$.
\begin{de}\label{def1}Any function $u\in  USC(\overline{\Om})$  (resp., $u\in  LSC(\overline{\Om})$) is called \emph{viscosity subsolution}
(resp., \emph{supersolution}) of
\begin{equation*}
\begin{cases}
 \di u+b(x)\cdot Du+c(x)u= g(x)  & \text{in} \quad\Om \\
 B(x,u,Du)= 0 & \text{on} \quad\partial\Om, \\
 \end{cases}
 \end{equation*}
 if the following conditions hold
\begin{itemize}
\item[(i)] For every $x_0\in \Om$, for all $\varphi\in C^2(\overline{\Om})$, such that $u-\varphi$ has a local maximum (resp., minimum)
at $x_0$ and $D\varphi(x_0)\neq0$, one has
$$\di\varphi(x_0)+b(x_0)\cdot D\varphi(x_0)+c(x_0)u(x_0)\geq \,(\text{resp., } \leq\,)\,\,g(x_0).$$ If
 $u\equiv k=$const. in a neighborhood of $x_0$, then
$$c(x_0) k\geq\, (\text{resp., } \leq\,)\,\,g(x_0).$$
\item[(ii)]For every $x_0\in \partial\Om$, for all $\varphi\in C^2(\overline{\Om})$, such that $u-\varphi$ has a local maximum (resp., minimum)
at $x_0$ and $D\varphi(x_0)\neq0$, one has
$$(-\di\varphi(x_0)-b(x_0)\cdot D\varphi(x_0)-c(x_0)u(x_0)+g(x_0))\wedge
B(x_0,u(x_0),D\varphi(x_0))\leq 0$$
(resp.,$$(-\di\varphi(x_0)-b(x_0)\cdot
D\varphi(x_0)-c(x_0)u(x_0)+g(x_0))\vee
B(x_0,u(x_0),D\varphi(x_0))\geq 0).$$ If $u\equiv k=$const. in a
neighborhood of $x_0$ in $\overline{\Om}$, then
$$(-c(x_0) k+g(x_0))\wedge B(x_0,k,0)\leq 0$$\quad (resp., $$(-c(x_0) k+g(x_0))\vee B(x_0,k,0)\geq 0).$$
\end{itemize}
\end{de}
It is possible to define sub and supersolutions of the
$\infty$-Laplace equation also using the semicontinous extensions
of the function $(p,X)\rightarrow \tr(\sigma(p)X)$ as done in
\cite{j} and \cite{jk}. In definition \ref{def1} it is remarkable
that nothing is required in the case $D\varphi(x_0)=0$ if $u$ is
not constant.

For a detailed presentation of the theory of viscosity solutions
and of the boun\-da\-ry conditions in the viscosity sense, we
refer the reader to e.g. \cite{cil}.

We call strong viscosity subsolutions (resp., supersolutions) the
viscosity subsolutions (resp., supersolutions) that satisfy
$B(x,u,Du)\leq $ (resp., $\geq$) 0 in the viscosity sense for all
$x\in\partial\Om$. If $\lam\rightarrow
B(x,r,p-\lam\overrightarrow{n})$ is non-increasing in $\lam\geq0$,
then classical subsolutions (resp., supersolutions) are strong
viscosity subsolutions (resp., supersolutions), see \cite{cil}
Proposition 7.2.

In the above definition the test functions can be substituted by
the elements of the semijets $\overline{J}^{2,+}u(x_0)$ when $u$
is a subsolution and $ \overline{J}^{2,-}u(x_0)$ when $u$ is a
supersolution, see \cite{cil}.
\section{Lipschitz continuity of viscosity solutions}
It is known that the $\infty$-harmonic functions, i.e., the
solution of $\di u=0$ are locally Lipschitz continuous, see e.g.
\cite{acj}. We now show the Lipschitz regularity in the whole
$\Oms$ of the solutions of the Neumann problem associated to the
$\infty$-Laplacian plus lower order terms.
\begin{thm}\label{regolaritathm}Assume that $\Om$ is a bounded domain of class $C^2$ and that $b$, $c$, $g$
are bounded in $\Om$. If $u\in C(\Oms)$ is a viscosity solution of
\begin{equation*}
\begin{cases}
 \di u+b(x)\cdot Du+c(x) u=  g(x) & \text{in} \quad\Om \\
 \frac{\partial u}{\partial\overrightarrow{n}}= 0 & \text{on} \quad\partial\Om, \\
 \end{cases}
 \end{equation*}then
$$ |u(x)-u(y)|\leq C_0|x-y|\quad \forall x,y\in\Oms,$$
where $C_0$ depends on
 $\Om,\,N,\,|b|_\infty,\,|c|_\infty,\,|g|_\infty,$ and
 $|u|_\infty$.
 \end{thm}
 The Theorem is an immediate consequence of the next lemma, the proof of which, though following the line of Proposition III.1
of \cite{il}, introduces new test functions that, in particular,
depend on the distance function $d(x)$.

The lemma will be used also in the proof of Theorem \ref{maxpneum}
in the next section.

\begin{lem}\label{regolarita}Assume the hypothesis of Theorem \ref{regolaritathm} and suppose that $g$ and $h$ are bounded functions.
Let $u\in USC(\Oms)$ be a viscosity subsolution of
\begin{equation*}
\begin{cases}
 \di u+b(x)\cdot Du+c(x) u=  g(x) & \text{in} \quad\Om \\
 \frac{\partial u}{\partial\overrightarrow{n}}= 0 & \text{on} \quad\partial\Om, \\
 \end{cases}
 \end{equation*}
and $v\in LSC(\Oms)$ a viscosity supersolution of
\begin{equation*}
\begin{cases}
 \di v+b(x)\cdot Dv+c(x) v=  h(x) & \text{in} \quad\Om \\
 \frac{\partial v}{\partial\overrightarrow{n}}= 0 & \text{on} \quad\partial\Om, \\
 \end{cases}
 \end{equation*}
with $u$ and $v$  bounded, or $v\geq 0$ and bounded. If
$m=\max_{\Oms}(u-v)\geq0$,
 then there exists $C_0>0$ such that
\begin{equation}\label{stimau-v1}u(x)-v(y)\leq m +C_0|x-y|\quad \forall
x,y\in\Oms,\end{equation}
 where $C_0$ depends on
 $\Om,\,N,\,|b|_\infty,\,|c|_\infty,\,|g|_\infty,\,|h|_\infty,\,|v|_\infty,\,m$ and
 $|u|_\infty$ or $\sup_{\Oms}u$.
 \end{lem}
 \dim
We set
$$\Phi(x)=MK|x|-M(K|x|)^2,$$and
$$\varphi (x,y)=m+e^{-L(d(x)+d(y))}\Phi(x-y),$$ where
$L$ is a fixed number greater than $2/(3r)$ with $r$ the radius in
the condition ($\Om2$) and where $K$ and $M$ are two positive
constants to be chosen later. If $K|x|\leq \frac{1}{4}$, then
\begin{equation}\label{phimagg}\Phi(x)\geq
\frac{3}{4}MK|x|.\end{equation} We define
$$\Delta _K:=\left\{(x,y)\in \R^N\times\R^N|\, |x-y|\leq\frac{1}{4K}\right\}.$$ We
fix $M$ such that
\begin{equation}\label{M}
\max_{\Oms^{\,2}}(u(x)-v(y))\leq
m+e^{-2Ld_0}\frac{M}{8},\end{equation} where
$d_0=\max_{x\in\Oms}d(x)$. To prove \eqref{stimau-v1} it is enough
to show that taking $K$ large enough, one has
$$u(x)-v(y)-\varphi(x,y)\leq 0\quad\text{for }(x,y)\in
\Delta_K\cap\Oms^2.$$ Suppose by contradiction that for each $K$
there is some point $(\xs,\ys)\in\Delta_K\cap \Oms^2$ such that
$$u(\xs)-v(\ys)-\varphi(\xs,\ys)=\max_{\Delta_K\cap
\Oms\,^2}(u(x)-v(y)-\varphi(x,y))>0.$$ Here we have dropped the
dependence of $\xs,\,\ys$ on $K$ for simplicity of notations.

Observe that if $v\geq0$, since from \eqref{phimagg} $\Phi(x-y)$
is non-negative in $\Delta_K$ and $m\geq 0$, one has $u(\xs)>0$.

Clearly $\xs\neq \ys$. Moreover the point $(\xs,\ys)$ belongs to
$\text{int}(\Delta_K)\cap\Oms^2$. Indeed, if $|x-y|=\frac{1}{4K}$,
by \eqref{M} and \eqref{phimagg} we have
\begin{equation*}
u(x)-v(y)\leq m +e^{-2Ld_0} \frac{M}{8}\leq
m+e^{-L(d(x)+d(y))}\frac{1}{2}MK|x-y|\leq\varphi(x,y).\end{equation*}
Since $\xs\neq \ys$ we can compute the derivatives of $\varphi$ at
$(\xs,\ys)$ obtaining
\begin{equation*}\begin{split}D_x\varphi(\xs,\ys)&=e^{-L(d(\xs)+d(\ys))}MK\Big\{-L|\xs-\ys|(1-K|\xs-\ys|)Dd(\xs)\\&+
(1-2K|\xs-\ys|)\frac{(\xs-\ys)}{|\xs-\ys|}\Big\},\end{split}\end{equation*}
\begin{equation*}\begin{split}D_y\varphi(\xs,\ys)&=e^{-L(d(\xs)+d(\ys))}MK\Big\{-L|\xs-\ys|(1-K|\xs-\ys|)Dd(\ys)\\&-
(1-2K|\xs-\ys|)\frac{(\xs-\ys)}{|\xs-\ys|}\Big\}.\end{split}\end{equation*}
Observe that for large $K$
\begin{equation}\label{dxphistima}0<e^{-L(d(\xs)+d(\ys))}MK\left(\frac{1}{2}-L|\xs-\ys|\right)\leq
|D_x\varphi(\xs,\ys)|,|D_y\varphi(\xs,\ys)|\leq 2MK.\end{equation}
Using \eqref{sferaest}, if $\xs\in\partial \Om$ we have
\begin{equation*}\begin{split}
&\langle D_x\varphi(\xs,\ys),\overrightarrow{n}(\xs)\rangle\\& =
e^{-Ld(\ys)}MK\Big\{L|\xs-\ys|(1-K|\xs-\ys|)+(1-2K|\xs-\ys|)\langle\frac{(\xs-\ys)}{|\xs-\ys|},\overrightarrow{n}(\xs)\rangle\Big\}\\&\geq
e^{-Ld(\ys)}MK\Big\{\frac{3}{4}L|\xs-\ys|-(1-2K|\xs-\ys|)\frac{|\xs-\ys|}{2r}\Big\}\\&\geq
\frac{1}{2}e^{-Ld(\ys)}MK|\xs-\ys|\left(\frac{3}{2}L-\frac{1}{r}\right)>0,
\end{split}\end{equation*}since
$\xs\neq\ys$ and $L>2/(3r)$. Similarly,
 if $\ys\in\partial\Om$
\begin{equation*}
\langle-D_y\varphi(\xs,\ys),\overrightarrow{n}(\ys)\rangle\leq
\frac{1}{2}e^{-Ld(\xs)}MK|\xs-\ys|\left(-\frac{3}{2}L+\frac{1}{r}\right)<0.\end{equation*}
In view of definition of sub and supersolution, we conclude that
$$\tr(\sigma(D_x\varphi(\xs,\ys))X)+b(\xs)\cdot D_x\varphi(\xs,\ys)+c(\xs)u(\xs)\geq g(\xs)\,\,\text{if }(D_x\varphi(\xs,\ys),X)\in
\overline{J}^{2,+}u(\xs),$$
$$\tr(\sigma(D_y\varphi(\xs,\ys))Y)-b(\ys)\cdot D_y\varphi(\xs,\ys)+c(\ys)v(\ys)\leq
h(\ys)\,\, \text{if }(-D_y\varphi(\xs,\ys),Y)\in
\overline{J}^{2,-}v(\ys).$$ Then the previous inequalities holds
for  any maximum point $(\xs,\ys)\in\Delta_K\cap\Oms\,^2$,
provided $K$ is large enough.

 Since $(\xs,\ys)\in\text{int}\Delta_K\cap\Oms\,^2$, it is a local
maximum of $u(x)-v(y)-\varphi(x,y)$ in $\Oms\,^2$. Applying
Theorem 3.2 in \cite{cil}, for every $\epsilon>0$ there exist
$X,Y\in \emph{S(N)}$ such that $ (D_x\varphi(\xs,\ys),X)\in
\overline{J}\,^{2,+}u(\xs),\,(-D_y\varphi(\xs,\ys),Y)\in
\overline{J}\,^{2,-}v(\ys)$ and
\begin{equation}\label{tm2ishii}\left(%
\begin{array}{cc}
  X & 0 \\
  0 & -Y \\
\end{array}%
\right)\leq D^2(\varphi(\xs,\ys))+\epsilon
(D^2(\varphi(\xs,\ys)))^2.
\end{equation}
Now we want to estimate the matrix on the right-hand side of the
last inequality.
\begin{equation*}\begin{split}D^2\varphi(\xs,\ys)&=\Phi(\xs-\ys)D^2(e^{-L(d(\xs)+d(\ys))})+D(e^{-L(d(\xs)+d(\ys))})\otimes
D(\Phi(\xs-\ys))\\&+D(\Phi(\xs-\ys))\otimes
D(e^{-L(d(\xs)+d(\ys))})+e^{-L(d(\xs)+d(\ys))}D^2(\Phi(\xs-\ys)).\end{split}\end{equation*}We
set $$A_1:=\Phi(\xs-\ys)D^2(e^{-L(d(\xs)+d(\ys))}),$$
$$A_2:=D(e^{-L(d(\xs)+d(\ys))})\otimes
D(\Phi(\xs-\ys))+D(\Phi(\xs-\ys))\otimes
D(e^{-L(d(\xs)+d(\ys))}),$$
$$A_3:=e^{-L(d(\xs)+d(\ys))}D^2(\Phi(\xs-\ys)).$$Observe that
\begin{equation}\label{a1}A_1\leq CK|\xs-\ys|\left(%
\begin{array}{cc}
  I & 0 \\
  0 & I \\
\end{array}%
\right).\end{equation} Here and henceforth C denotes various
positive constants independent of $K$.

For $A_2$ we have the following estimate
\begin{equation}\label{a2}A_2\leq
CK\left(%
\begin{array}{cc}
  I & 0 \\
  0 & I \\
\end{array}%
\right)
+CK\left(%
\begin{array}{cc}
  I & -I \\
  -I & I \\
\end{array}%
\right).\end{equation} Indeed for $\xi,\,\eta\in\R^N$ we compute
\begin{equation*}\begin{split}\langle
A_2(\xi,\eta),(\xi,\eta)\rangle&=2Le^{-L(d(\xs)+d(\ys))}\{\langle
Dd(\xs)\otimes D\Phi(\xs-\ys)(\eta-\xi),\xi\rangle\\&+\langle
Dd(\ys)\otimes D\Phi(\xs-\ys)(\eta-\xi),\eta\rangle\}\leq
CK(|\xi|+|\eta|)|\eta-\xi|\\&\leq
CK(|\xi|^2+|\eta|^2)+CK|\eta-\xi|^2.\end{split}\end{equation*}

Now we consider $A_3$. The matrix $D^2(\Phi(\xs-\ys))$ has the
form
$$D^2(\Phi(\xs-\ys))=\left(%
\begin{array}{cc}
  D^2\Phi(\xs-\ys) & - D^2\Phi(\xs-\ys) \\
  - D^2\Phi(\xs-\ys) &  D^2\Phi(\xs-\ys) \\
\end{array}%
\right),$$and the Hessian matrix of $\Phi(x)$ is
\begin{equation}\label{hessianphi}D^2\Phi(x)=\frac{MK}{|x|}\left(I-\frac{x\otimes
x}{|x|^2}\right)-2MK^2I.\end{equation} If we choose
\begin{equation}\label{epsilon}\epsilon=\frac{|\xs-\ys|}{2MKe^{-L(d(\xs)+d(\ys))}},\end{equation}
then we have the following estimates
$$\epsilon A_1^2\leq
CK|\xs-\ys|^3I_{2N},\quad \epsilon A_2^2\leq CK|\xs-\ys|I_{2N},$$
\begin{equation}\label{aprodotti}\begin{split} \epsilon (A_1A_2+A_2A_1)\leq
CK|\xs-\ys|^2I_{2N},\end{split}\end{equation}
$$\epsilon (A_1A_3+A_3A_1)\leq CK|\xs-\ys|I_{2N},\quad \epsilon
(A_2A_3+A_3A_2)\leq CKI_{2N},$$
 where $I_{2N}:=\left(%
\begin{array}{cc}
  I & 0 \\
  0 & I \\
\end{array}%
\right)$. Then using \eqref{a1}, \eqref{a2},
\eqref{aprodotti} and observing that $$(D^2(\Phi(\xs-\ys)))^2=\left(%
\begin{array}{cc}
  2(D^2\Phi(\xs-\ys))^2 & - 2(D^2\Phi(\xs-\ys))^2 \\
  - 2(D^2\Phi(\xs-\ys))^2 &  2(D^2\Phi(\xs-\ys))^2 \\
\end{array}%
\right),$$from \eqref{tm2ishii} we conclude that
$$\left(%
\begin{array}{cc}
  X & 0 \\
  0 & -Y \\
\end{array}%
\right)\leq O(K)\left(%
\begin{array}{cc}
  I & 0 \\
  0 & I \\
\end{array}%
\right)+\left(%
\begin{array}{cc}
  B & -B \\
  -B & B \\
\end{array}%
\right),$$ where
\begin{equation*}B=CKI+e^{-L(d(\xs)+d(\ys))}\left[D^2\Phi(\xs-\ys)+\frac{|\xs-\ys|}{MK}
(D^2\Phi(\xs-\ys))^2\right].\end{equation*}The last inequality can
be rewritten as follows
$$\left(%
\begin{array}{cc}
  \widetilde{X} & 0 \\
  0 & -\widetilde{Y} \\
\end{array}%
\right)\leq\left(%
\begin{array}{cc}
  B & -B \\
  -B & B \\
\end{array}%
\right),$$ with $\widetilde{X}=X-O(K)I$ and
$\widetilde{Y}=Y+O(K)I.$ Multiplying on the left the previous
inequality by the non-negative symmetric matrix
$$\left(%
\begin{array}{cc}\sigma(D_x\varphi(\xs,\ys)) & 0\\ 0 &
\sigma(D_y\varphi(\xs,\ys))\end{array}%
\right),$$ and taking traces we get
\begin{equation}\label{fu-fv}\tr(\sigma(D_x\varphi(\xs,\ys))\widetilde{X})-\tr(\sigma(D_y\varphi(\xs,\ys))\widetilde{Y})\leq
\tr(\sigma(D_x\varphi(\xs,\ys))B)+\tr(\sigma(D_y\varphi(\xs,\ys))B).\end{equation}

We want to get a good estimate for the matrix on the right-hand
side above. For that aim let
$$0\leq P:=\frac{(\xs-\ys)\otimes (\xs-\ys)}{|\xs-\ys|^2}\leq I,$$
 and let us compute tr($PB$). From \eqref{hessianphi}, since
the matrix $(1/|x|^2)x\otimes x$ is idempotent, we get
$$(D^2\Phi(x))^2=\frac{M^2K^2}{|x|^2}(1-4K|x|)\left(I-\frac{x\otimes
x}{|x|^2}\right)+4M^2K^4I.$$ Then, using that $\text{tr}P=1$ and
$4K|\xs-\ys|\leq1$, we have
\begin{equation*}\begin{split}\text{tr}(PB)&=CK+e^{-L(d(\xs)+d(\ys))}(-2MK^2+4MK^3|\xs-\ys|)
\\&\leq CK-e^{-L(d(\xs)+d(\ys))}MK^2\leq -CK^2, \end{split}\end{equation*}for
large $K$. The vector $D_x\varphi(\xs,\ys)$ can be written in the
following way
\begin{equation*}\begin{split}D_x\varphi(\xs,\ys)&=e^{-L(d(\xs)+d(\ys))}MK(v_1+v_2),\end{split}\end{equation*}where
$$v_1=-L|\xs-\ys|(1-K|\xs-\ys|)Dd(\xs),\quad
v_2=(1-2K|\xs-\ys|)\frac{(\xs-\ys)}{|\xs-\ys|},$$and so
$$\sigma(D_x\varphi(\xs,\ys))=\frac{v_1\otimes
v_1}{|v_1+v_2|^2}+\frac{v_1\otimes v_2+v_2\otimes
v_1}{|v_1+v_2|^2}+\frac{v_2\otimes v_2}{|v_1+v_2|^2}.$$  Since
$K|\xs-\ys|\leq\frac{1}{4}$, for large $K$ we have
$$\frac{1}{4}=\frac{1}{2}-\frac{1}{4}\leq|v_2|-|v_1|\leq
|v_1+v_2|\leq |v_1|+|v_2|\leq 2,$$ and
$$\|B\|\leq \frac{CK}{|\xs-\ys|}.$$Then
$$\left|\tr\left(\frac{v_1\otimes
v_1}{|v_1+v_2|^2}B\right)\right|\leq C|\xs-\ys|^2\|B\|\leq
CK|\xs-\ys|,$$
$$\left|\tr\left(\frac{v_1\otimes v_2+v_2\otimes
v_1}{|v_1+v_2|^2}B\right)\right|\leq C|\xs-\ys|\|B\|\leq CK$$and
$$\tr\left(\frac{v_2\otimes
v_2}{|v_1+v_2|^2}B\right)=\frac{1}{|v_1+v_2|^2}\tr(PB)\leq
 -CK^2.$$ In conclusion
$$\tr(\sigma(D_x\varphi(\xs,\ys)B))\leq O(K)-CK^2.$$ The same
estimate holds for $\tr(\sigma(D_y\varphi(\xs,\ys))B).$ Hence,
from \eqref{fu-fv} we conclude that
$$\tr(\sigma(D_x\varphi(\xs,\ys)\widetilde{X})-\tr(\sigma(D_y\varphi(\xs,\ys)\widetilde{Y})\leq
O(K)-CK^2.$$ Now, using the previous estimate, the definition of
$\widetilde{X}$ and $\widetilde{Y}$ and the fact that $u$ and $v$
are respectively sub and supersolution we
compute\begin{equation*}\begin{split}g(\xs)-c(\xs) u(\xs)&\leq
\text{tr}(\sigma(D_x\varphi)X)+b(\xs)\cdot D_x\varphi\\&\leq
\text{tr}(\sigma(D_x\varphi)\widetilde{X}) + O(K)+b(\xs)\cdot
D_x\varphi\\& \leq \text{tr}(\sigma(D_y\varphi){Y})
+O(K)-CK^2+b(\xs)\cdot D_x\varphi\\&\leq b(\ys)\cdot D_y\varphi
-c(\ys) v(\ys)+h(\ys) + O(K)-CK^2+b(\xs)\cdot
D_x\varphi.\end{split}\end{equation*} From this inequalities,
using \eqref{dxphistima} we get
\begin{equation*}
g(\xs)-h(\ys)-c(\xs) u(\xs)+c(\ys)v(\ys)\leq
O(K)-CK^2.\end{equation*} If both $u$ and $v$ are bounded, then
the
 member on the left-hand side of the last inequality is bounded from below by
$-|g|_\infty-|h|_\infty-|c|_\infty(|u|_\infty+|v|_\infty)$.
Otherwise, if $v$ is non-negative and bounded, then $u(\xs)\geq 0$
and that quantity is greater than
$-|g|_\infty-|h|_\infty-|c|_\infty (\sup u+ |v|_\infty)$. On the
other hand, the  member on the right-hand side goes to $-\infty$
as $K\rightarrow+\infty$, hence taking $K$ large enough we obtain
a contradiction and this concludes the proof.\finedim
\begin{rem}\label{stimadirich}{\em If   $u$ is a subsolution of $\di u+b(x)\cdot Du+c(x)u=g$,
$v$ is a supersolution of $\di v+b(x)\cdot Dv+c(x)v=h$ in $\Om$,
$u\leq v$ on $\partial\Om$ and $m>0$ then the estimate
\eqref{stimau-v1} still holds for any $x,y\in\Om$. To prove this
define $\varphi=m+MK|x|-M(K|x|)^2$ and follow the proof of Lemma
\ref{regolarita}.}
\end{rem}

Since the Lipschitz estimate  depends only on the bounds of the
solution of $g$ and on the structural constants, an immediate
consequence of Theorem \ref{regolaritathm} is the following
compactness criterion that will be useful in the next sections.
\begin{cor}\label{corcomp}Assume the hypothesis of Theorem \ref{regolaritathm} on $\Om$, $F$ and $b$. Suppose that $(g_n)_{n}$ is a
sequence of continuous and uniformly bounded functions and
$(u_n)_n$ is a sequence of uniformly bounded viscosity solutions
of
\begin{equation*}
\begin{cases}
 \di u_n+b(x)\cdot Du_n=  g_n(x) & \text{in} \quad\Om \\
 \frac{\partial u_n}{\partial\overrightarrow{n}}= 0 & \text{on} \quad\partial\Om. \\
 \end{cases}
 \end{equation*}Then the sequence $(u_n)_n$ is relatively compact in
 $C(\Oms)$.
 \end{cor}

\section{The Maximum Principle and the principal eigenvalues}
\noindent We say that the operator $ \di+b(x)\cdot D+c(x)$ with
the Neumann boundary condition satisfies the maximum principle if
whenever $u\in USC(\Oms)$ is a viscosity subsolution of
\begin{equation*}
\begin{cases}
 \di u+b(x)\cdot Du+c(x)u= 0  & \text{in} \quad\Om \\
 \frac{\partial u}{\partial\overrightarrow{n}}= 0 & \text{on} \quad\partial\Om, \\
 \end{cases}
 \end{equation*}then $u\leq 0$ on $\Oms$.

We first prove that the maximum principle holds under the
classical assumption $c\leq 0$, also for domain which are not of
class $C^2$ and with more general boundary conditions. Then we
show that the operator $\di +b(x)\cdot D+c(x)+\lam $ with the
Neumann boundary condition satisfies the maximum principle for any
$\lam<\lams$. This is the best result that one can expect, indeed,
as we will see, $\lams$ admits a positive eigenfunction which
provides a counterexample to the maximum principle for
$\lam\geq\lams$.

Finally, we give an example of class of functions $c(x)$ which
change sign in $\Om$ and such that the associated principal
eigenvalue $\lams$ is positive.

 \subsection{The case $c(x)\leq 0$}
In this subsection we assume that $\Om$ is of class $C^1$ and
satisfies the interior sphere condition ($\Om$1). We need the
comparison principle between sub and supersolutions of the
Dirichlet problem when $c<0$ in $\Om$. This result is known for
the operator $\di u+b(x)\cdot Du+c(x)u$ when $b$ is Lipschitz
continuous or $b$ satisfies $\langle b(x)-b(y),x-y\rangle\leq0$,
see e.g. \cite{cil}. Actually, we can remove these conditions.
\begin{thm}\label{dirichletcomp} Let $\Om$ be bounded. Assume that $b$, $c$ and $g$ are continuous and bounded in
$\Om$ and $c<0$ on $\Oms$. If $u\in USC(\Oms)$ and $v\in
LSC(\Oms)$ are respectively sub and supersolution of
$$ \di u+b(x)\cdot Du+c(x)u= g(x) \quad\text{in } \Om,$$ and $u\leq
v$ on $\partial\Om$ then $u\leq v$ in $\Om$.
\end{thm}
For convenience of the reader the proof of the theorem will be
sketched at the end of the next subsection.

The previous comparison result allows us to establish  the strong
minimum and maximum principles, for sub and supersolutions of the
Neumann problem even with the following more general boundary
condition $$f(x,u)+\frac{\partial u}{\partial\overrightarrow{n}}=0
\quad x\in\partial\Om,$$for some
$f:\partial\Om\times\R\rightarrow\R$.
\begin{prop}\label{p1}Let $\Om$ be a $C^1$ domain satisfying ($\Om$1). Suppose that
$b$ and $c$ are bounded and continuous in $\Om$ and that
$f(x,0)\leq 0$ for all $x\in\partial\Om$. If $v\in
LSC(\overline{\Om})$ is a non-negative viscosity supersolution of
\begin{equation}\label{moregencond}
\begin{cases}
\di v+b(x)\cdot Dv+c(x)v= 0  & \text{in} \quad\Om \\
 f(x,v)+\frac{\partial v}{\partial\overrightarrow{n}}= 0 & \text{on} \quad\partial\Om, \\
 \end{cases}
 \end{equation}
 then either $v\equiv0$ or $v>0$ on $\overline{\Om}$.
\end{prop}
\dim Since $v$ is non-negative, it is supersolution in $\Om$ of
the equation
\begin{equation}\label{p1e2}
\di v+b(x)\cdot Dv-|c|_\infty v= 0.
\end{equation}
Without loss of generality we can assume $|c|_\infty>0$. Suppose
by contradiction that $v\not\equiv 0$ vanishes somewhere in $\Om$.
Then we can find $x_1,x_0\in\Om$ and $R>0$ such that
$B(x_1,\frac{3}{2}R)\subset\Om$, $v>0$ in $B(x_1,R)$,
$|x_1-x_0|=R$ and $v(x_0)=0$. Let us construct a subsolution of
\eqref{p1e2} in the annulus $\frac{R}{2}< |x-x_1|=r<
\frac{3}{2}R.$

Let us consider the function $\phi(x)=e^{-kr}-e^{-kR}$, where $k$
is a positive constant to be determined. It easy to see that for
radial functions $g(x)=\varphi(r)$, $\di g(x)=\varphi^{''}(r)$.
Then
\begin{equation*}\begin{split} \di \phi+b(x)\cdot D\phi-|c|_{\infty}\phi&
=k^2e^{-kr}-ke^{-kr}b(x)\cdot
\dfrac{(x-y)}{r}-|c|_{\infty}(e^{-kr}-e^{-kR})\\
&\geq e^{-kr}\left (k^{2}-|b|_\infty k
-|c|_\infty\right).\end{split}\end{equation*} Take $k$ such that
$$k^{2}-|b|_\infty k
-|c|_\infty>0,$$ then $\phi$ is a strict subsolution of the
equation \eqref{p1e2}. Now choose $m>0$ such that
$$m(e^{-k\frac{R}{2}}-e^{-kR})=v_1:=\text
{inf}_{|x-x_1|=\frac{R}{2}}v(x)>0,$$ and define
$w(x)=m(e^{-kr}-e^{-kR})$. By homogeneity $w$ is still a
subsolution of \eqref{p1e2} in the annulus
$\frac{R}{2}<|x-x_1|<\frac{3}{2}R$, moreover $w=v_1\leq v $ if
$|x-x_1|=\frac{R}{2}$ and $w< 0 \leq v$ if $|x-x_1|=\frac{3}{2}R$.
Then by the comparison principle, Theorem \ref{dirichletcomp},
$w\leq v$ in the entire annulus.

Since $v(x_0)=w(x_0)=0$, $w$ is a test function for $v$ at $x_0$
with $Dw(x_0)\neq 0$. But $$\di w(x_0)+b(x_0)\cdot
Dw(x_0)-|c|_\infty v(x_0)> 0,$$ and this contradicts the
definition of $v$. Then $v>0$ in $\Om$.

Now suppose by contradiction that $x_0$ is some point in $\partial
\Om$ on which $v(x_0)=0$. The interior sphere condition ($\Om$1)
implies that there exist $R>0$ and $y\in \Om$ such that the ball
centered in $y$ and of radius $R$, $B(y,R)$, is contained in $\Om$
and $x_0\in
\partial B(y,R)$. Fixed $0<\rho<R$, as before the function
$w(x)=m(e^{-kr}-e^{-kR})$ is a strict subsolution of  \eqref{p1e2}
in the annulus $\rho<|x-y|=r<R$, where $m$ is such that
$m(e^{-k\rho}-e^{-kR})=v_1:=\text {inf}_{|x-y|=\rho }v(x)>0$.
Since $w\leq v$ on the boundary of the annulus then again by the
comparison principle, Theorem \ref{dirichletcomp}, $w\leq v$ in
the entire annulus.

 Now let $\delta$ be a
positive number smaller than $R-\rho$. In
$B(x_0,\delta)\cap\overline{\Om}$ still $w\leq v$, since for
$|x-y|>R$, $w< 0 \leq v$; moreover $w(x_0)=v(x_0)=0$. Then $w$ is
a test function for $v$ at $x_0$. But
$$\di w(x_0)+b(x_0)\cdot
Dw(x_0)-|c|_{\infty}v(x_0)> 0,$$ and
$$f(x_0,v(x_0))+ \dfrac{\partial w}{\partial \overrightarrow{n}}(x_0)=f(x_0,0)-kme^{-kR}<0.$$ This contradicts the definition of $v$.
Finally $v$ cannot be zero on $\overline{\Om}$.\finedim

Similarly we can prove
\begin{prop}\label{p2}Let $\Om$ be a $C^1$ domain satisfying ($\Om$1). Assume that
$b$ and $c$ are bounded and continuous in $\Om$ and that
$f(x,0)\geq 0$ for all $x\in\partial\Om$. If $u\in
USC(\overline{\Om})$ is a non-positive viscosity subsolution of
\eqref{moregencond} then either $u\equiv0$ or $u<0$ on
$\overline{\Om}$.
\end{prop}
For $x\in\partial\Om$, let us introduce $S(x)$, the symmetric
operator corresponding to the second fundamental form of $\partial
\Om$ in $x$ oriented with the exterior normal to $\Om$.

\begin{thm}[Maximum Principle for $c\leq 0$]\label{macpcless0}Assume the hypothesis of Proposition \ref{p2}. In addition suppose that
$\Om$ is bounded, $c\leq 0$, $c\not\equiv 0$ and $r\rightarrow
f(x,r)$ is non-decreasing on $\R$. If $u\in USC(\overline{\Om})$
is a viscosity subsolution of \eqref{moregencond} then $u\leq 0$
on $\Oms$. The same conclusion holds also if $c\equiv 0$ in the
following two cases
\begin{itemize}
\item[(i)] $\Om$ is a $C^2$ domain and for any $r>0$ there exists $\xs\in \partial\Om$
 such that $f(\xs,r)>0$, $S(\xs)\leq 0$ and $\langle b(\xs), \overrightarrow{n}(\xs)\rangle>0$;
\item[(ii)]$\max_{x\in\partial\Om}f(x,r)>0$ for any $r>0$ and $u$ is a strong subsolution.
\end{itemize}
\end{thm}
\dim Let $u$ be a subsolution of \eqref{moregencond} and
$c\not\equiv 0$. First let us suppose $u\equiv k=$const. By
definition
$$c(x)k\geq 0\qquad \text{in }\Om,$$ which implies $k\leq
0$.

Now we assume that $u$ is not a constant. We argue by
contradiction; suppose that $\max_{\Oms}u=u(x_0)>0$, for some
$x_0\in\Oms$. Define $\widetilde{u}(x):=u(x)-u(x_0)$. Since $c\leq
0$ and $f$ is non-decreasing, $\widetilde{u}$ is a non-positive
subsolution of \eqref{moregencond}. Then, from Proposition
\ref{p2}, either $u\equiv u(x_0)$ or $u<u(x_0)$ on
$\overline{\Om}$. In both cases we get a contradiction.

Let us turn to the case $c\equiv 0$. We have to prove that $u$
cannot be a positive constant. Suppose by contradiction that
$u\equiv k$. Suppose that $\Om$ is a $C^2$ domain and let $\xs\in
\partial \Om$ be such that
$S(\xs)\leq 0$, $\langle b(\overline{x}),
\overrightarrow{n}(\overline{x})\rangle>0$ and
$f(\overline{x},k)>0$. In general, if $\phi$ is a $C^2$ function,
$\overline{x}\in
\partial \Om$ and $S(\xs)\leq 0$, then
$(D\phi(\overline{x})-\lambda \overrightarrow{n}(\xs),
D^2\phi(\overline{x}))\in J^{2,+}\phi(\overline{x})$, for $\lambda
\geq0$ (see \cite{cil} Remark 2.7). Hence $(-\lambda
\overrightarrow{n}(\overline{x}),0)\in J^{2,+}u(\overline{x})$.
But $$f(\overline{x},k)-\lam\langle
\overrightarrow{n}(\overline{x}),\overrightarrow{n}(\overline{x})\rangle=
f(\overline{x},k)-\lambda>0,$$ for $\lambda>0$ small enough, and
$$-\lambda\langle b(\overline{x}), \overrightarrow{n}(\overline{x})\rangle<0.$$
This contradicts the definition of $u$.

Finally if $u$ is a strong subsolution, $u\equiv k>0$ and
$f(\overline{x},k)>0$ for some $\overline{x}\in\partial\Om$, then
the boundary condition is not satisfied at $\xs$ for $p=0$.
\finedim
\begin{rem}{\em Under the same assumptions of Theorem
\ref{macpcless0}, but now with $f$ sa\-ti\-sfying $f(x,0)\leq 0$
for all $x\in\partial\Om$ and with $f(\xs,r)<0$ for  $r<0$ in (i)
and $\min_{x\in\partial\Om}f(x,r)<0$ for $r<0$ in (ii), using
Proposition \ref{p1} we can prove the minimum principle, i.e., if
$u\in LSC(\Oms)$ is a viscosity supersolution of
\eqref{moregencond} then $u\geq 0$ on $\Oms$.}\end{rem}
\begin{rem}\emph{$C^2$ convex sets satisfy the condition $S\leq 0$ in
every point of the boundary.}\end{rem}
\begin{rem}{\em If $c\equiv 0$ and $f\equiv 0$ a counterexample to the maximum
principle is given by the positive constants.}\end{rem}

\subsection{The threshold for the Maximum Principle}

\noindent In this subsection and in the rest of the paper we
always assume that $\Om$ is bounded and of class $C^2$ and that
$b$ and $c$ are continuous on $\Oms$.

\begin{thm}[Maximum Principle for $\lam<\lams$]\label{maxpneum} Let $\lam <\lams$ and let $u\in
USC(\Oms)$ be a viscosity subsolution of
\begin{equation}\label{equazprincmax}
\begin{cases}
 \di u+b(x)\cdot Du+(c(x)+\lam)u = 0  & \text{in} \quad\Om \\
 \dfrac{\partial u}{\partial \overrightarrow{n}}= 0 & \text{on} \quad\partial\Om, \\
 \end{cases}
 \end{equation}then $u\leq 0$ on $\Oms$.
 \end{thm}
\begin{cor}The quantity $\lams$ is
finite.
\end{cor}
\dim It suffices to observe that $\lams\leq |c|_{\infty}$, since
when the zero order coefficient is $c(x)+|c|_\infty$ the maximum
principle does not hold. A counterexample is given by the positive
constants. \finedim

In the proof of Theorem \ref{maxpneum}  we need the following
result which is an adaptation of Lemma 1 of \cite{bd} for
supersolutions of the Neumann boundary value problem.
\begin{lem}\label{lemxdivy}Let $v\in LSC(\Oms)$ be a viscosity
supersolution of
\begin{equation*}
\begin{cases} \di v+b(x)\cdot Dv-\beta(v(x))=
g(x) & \text{in} \quad\Om \\
 \dfrac{\partial v}{\partial \overrightarrow{n}}= 0 & \text{on} \quad\partial\Om, \\
 \end{cases}
 \end{equation*} for some functions $g,\beta\in USC(\Oms)$. Suppose that $\xs\in \Oms$
is a strict local minimum of $v(x)+C|x-\xs|^qe^{-kd(x)}$,
$k>\frac{q}{2r}$, where $r$ is the radius in the condition
($\Om2$) and $q>2$. Moreover suppose that $v$ is not locally
constant around $\xs$. Then
$$-\beta(v(\xs))\leq g(\xs).$$\end{lem}

\begin{rem}\emph{Similarly, if $\beta,$ $g\in LSC(\Oms)$, $u\in USC(\Oms)$ is a
supersolution, $\overline{x}$ is a strict local maximum of
$u(x)-C|x-\xs|^qe^{-kd(x)}$, $k>\frac{q}{2r}$, $q>2$ and $u$ is
not locally constant around $\xs$, it can be proved that
$$-\beta(u(\xs))\geq g(\xs).$$}\end{rem}

 \dims \textbf{of Theorem \ref{maxpneum}.}
 Let $\tau\in ]\lam,\lams[$, then by definition there exists
 $v>0$ on $\Oms$ bounded viscosity supersolution of
 \begin{equation}\label{soprasol1}
 \begin{cases}
 \di v+b(x)\cdot Dv+(c(x)+\tau)v =0  & \text{in} \quad\Om \\
 \dfrac{\partial v}{\partial \overrightarrow{n}} = 0 & \text{on} \quad\partial\Om. \\
 \end{cases}
 \end{equation}

 We argue by contradiction and suppose that $u$ has a positive maximum in
 $\Oms$.
As in \cite{bd}, we define $\gamma':=\text{sup}_{\Oms}(u/v)>0$ and
$w=\gamma v$, with $\gamma\in (0,\gamma')$ to be determined. By
homogeneity, $w$ is still a supersolution of \eqref{soprasol1}.
Let $\overline{y}\in\Oms$ be such that
$u(\overline{y})/v(\overline{y})=\gamma'$. Since
$u(\overline{y})-w(\overline{y})=(\gamma'
-\gamma)v(\overline{y})>0$, the supremum of $u-w$ is strictly
positive, then by upper semicontinuity there exists
$\overline{x}\in\Oms$
 such that
$$ u(\xs)-w(\xs)=\max_{\Oms}(u-w)=m>0.$$
Clearly $u(\xs)>w(\xs)>0,$ moreover $u(\xs)\leq
\gamma'v(\xs)=\frac{\gamma'}{\gamma}w(\xs),$ from which
\begin{equation}\label{w(x)}
w(\xs)\geq \frac{\gamma}{\gamma'}u(\xs).
\end{equation}

Fix $q>2$ and $k>q/(2r)$, where $r$ is the radius in the condition
($\Om2$), and define for $j\in\N$ the functions $\phi\in
C^2(\Oms\times\Oms)$ and $\psi\in USC(\Oms\times\Oms)$ by
$$\phi(x,y)=\frac{j}{q}|x-y|^qe^{-k(d(x)+d(y))},\quad\psi(x,y)=u(x)-w(y)-\phi(x,y).$$
Let $(x_j,y_j)\in \Oms\times\Oms$ be a maximum point of $\psi$,
then $m=\psi(\xs,\xs)\leq u(x_j)-w(y_j)-\phi(x_j,y_j)$, from which
\begin{equation}\label{dis1}
\frac{j}{q}|x_j-y_j|^q\leq
(u(x_j)-w(y_j)-m)e^{k(d(x_j)+d(y_j))}\leq C,
\end{equation} where $C$ is independent
of $j$. The last relation implies that, up to subsequence, $x_j$
and $y_j$ converge to some $\zs\in\Oms$ as $j\rightarrow +\infty$.
Classical arguments show that
$$\lim_{j\rightarrow+\infty}\frac{j}{q}|x_j-y_j|^q=0,\quad\lim_{j\rightarrow+\infty}u(x_j)=u(\zs),\quad
\lim_{j\rightarrow+\infty}w(y_j)=w(\zs),$$and $$u(\zs)-w(\zs)=m.$$

\textbf{Claim 1} \emph{For $j$ large enough, there exist $x_j$ and
$y_j$ such that $(x_j,y_j)$ is a maximum point of $\psi$ and
$x_j\neq y_j$.}

Indeed if $x_j=y_j$ we have
$$\psi(x_j,x)=u(x_j)-w(x)-\frac{j}{q}|x-x_j|^qe^{-k(d(x_j)+d(x))}\leq\psi(x_j,x_j)=u(x_j)-w(x_j),$$and
$$\psi(x,x_j)=u(x)-w(x_j)-\frac{j}{q}|x-x_j|^qe^{-k(d(x)+d(x_j))}\leq\psi(x_j,x_j)=u(x_j)-w(x_j).$$ Then $x_j$ is a
minimum point for
$$W(x):=w(x)+\frac{j}{q}e^{-kd(x_j)}|x-x_j|^qe^{-kd(x)},$$and a
maximum point for
$$U(x):=u(x)-\frac{j}{q}e^{-kd(x_j)}|x-x_j|^qe^{-kd(x)}.$$
We first exclude that $x_j$ is both a strict local minimum and a
strict local maximum. Indeed in that case, if $u$ and $w$ are not
locally constant around $x_j$, by Lemma \ref{lemxdivy}
\begin{equation*}(c(x_j)+\tau)w(x_j)\leq (c(x_j)+\lam )u(x_j).\end{equation*}
The same result holds if $u$ or $w$ are locally constant by
definition of sub and supersolution. The last inequality leads to
a contradiction, as we will see at the end of the proof. Hence
$x_j$ cannot be both a strict local minimum and a strict local
maximum. In the first case there exist $\delta
>0$ and $ R>\delta$ such that
\begin{equation*}\begin{split}w(x_j)&=\min_{\delta\leq|x-x_j|\leq R\atop
x\in\Oms}\left(w(x)+\frac{j}{q}|x-x_j|^qe^{-k(d(x_j)+d(x))}\right)\\&=w(y_j)+\frac{j}{q}|y_j-x_j|^qe^{-k(d(x_j)+d(y_j))},\end{split}\end{equation*}
for some $y_j\neq x_j$, so that $(x_j,y_j)$ is still a maximum
point for $\psi$. In the other case, similarly, one can replace
$x_j$ by a point $y_j\neq x_j$ such that $(y_j,x_j)$ is a maximum
for $\psi$. This concludes the Claim 1.

Now computing the derivatives of $\phi$ we get
$$D_x\phi(x,y)=j|x-y|^{q-2}e^{-k(d(x)+d(y))}(x-y)-k\frac{j}{q}|x-y|^{q}e^{-k(d(x)+d(y))}Dd(x),$$
and
$$D_y\phi(x,y)=-j|x-y|^{q-2}e^{-k(d(x)+d(y))}(x-y)-k\frac{j}{q}|x-y|^{q}e^{-k(d(x)+d(y))} Dd(y).$$
Denote $p_j:=D_x\phi(x_j,y_j)$ and $r_j:=-D_y\phi(x_j,y_j).$ Since
$x_j\neq y_j$, $p_j$ and $r_j$ are different from 0 for $j$ large
enough. Indeed
\begin{equation}\label{graddiff0}0<\frac{j}{2}|x_j-y_j|^{q-1}e^{-2kd_0}\leq|p_j|,|r_j|\leq
2j|x_j-y_j|^{q-1},\end{equation} for large $j$, where
$d_0=\max_{\Oms}d(x)$. Using \eqref{sferaest}, if $x_j\in\partial
\Om$ then
\begin{equation*}
\langle p_j,\overrightarrow{n}(x_j)\rangle\geq j
|x_j-y_j|^{q}e^{-kd(y_j)}\left(-\frac{1}{2r}+\frac{k}{q}\right)>0,\end{equation*}
and if $y_j\in\partial\Om$ then
\begin{equation*}
\langle r_j,\overrightarrow{n}(y_j)\rangle\leq j
|x_j-y_j|^{q}e^{-kd(x_j)}\left(\frac{1}{2r}-\frac{k}{q}\right)<0,\end{equation*}
since $k>q/(2r)$ and $x_j\neq y_j$. In view of definition of sub
and supersolution we conclude that
$$\tr(\sigma(p_j)X)+b(x_j)\cdot p_j +(c(x_j)+\lam)
u(x_j)\geq 0\quad \text{if }(p_j,X)\in \overline{J}^{2,+}u(x_j),$$
$$\tr(\sigma(r_j)Y)+b(y_j)\cdot r_j+(c(y_j)+\tau) w(y_j)\leq 0\quad\text{if
}(r_j,Y)\in \overline{J}^{2,-}w(y_j).$$

Applying Theorem 3.2 of \cite{cil} for any $\epsilon>0$ there
exist $X_j, Y_j\in \emph{S(N)}$ such that $(p_j,X_j)\in
\overline{J}^{2,+}u(x_j)$, $(r_j,Y_j)\in \overline{J}^{2,-}w(y_j)$
and

 \begin{equation}\label{tm1lemis}-\left(\frac{1}{\epsilon}+\|D^2\phi(x_j,y_j)\|\right)\left(%
\begin{array}{cc}
  I & 0 \\
  0 & I \\
\end{array}%
\right)\leq \left(%
\begin{array}{cc}
  X_j & 0 \\
  0 & -Y_j \\
\end{array}%
\right)\leq D^2\phi(x_j,y_j)+\epsilon (D^2\phi(x_j,y_j))^2.
\end{equation}

\textbf{Claim 2} \emph{$X_j$ and $Y_j$ satisfy
\begin{equation}\label{claimmatrix} \left(%
\begin{array}{cc}
  X_j-\widetilde{X_j} & 0 \\
  0 & -Y_j+\widetilde{Y_j} \\
\end{array}%
\right)\leq \zeta_j\left(%
\begin{array}{cc}
  I & -I \\
  -I & I \\
\end{array}%
\right),\end{equation} where $\zeta_j=Cj|x_j-y_j|^{q-2}$, for some
positive constant $C$ independent of $j$ and  some matrices
$\widetilde{X_j},\,\widetilde{Y_j}=O(j|x_j-y_j|^{q}).$}

To prove the claim we need to estimate $D^2\phi(x_j,y_j)$.
\begin{equation*}\begin{split}D^2\phi(x_j,y_j)&=\frac{j}{q}|x_j-y_j|^qD^2(e^{-k(d(x_j)+d(y_j))})+
D(e^{-k(d(x_j)+d(y_j))})\otimes
\frac{j}{q}D(|x_j-y_j|^q)\\&+\frac{j}{q}D(|x_j-y_j|^q)\otimes
D(e^{-k(d(x_j)+d(y_j))})+e^{-k(d(x_j)+d(y_j))}\frac{j}{q}D^2(|x_j-y_j|^q).\end{split}\end{equation*}
We denote
$$A_1:=\frac{j}{q}|x_j-y_j|^qD^2(e^{-k(d(x_j)+d(y_j))}),$$
$$A_2:=De^{-k(d(x_j)+d(y_j))}\otimes
\frac{j}{q}D(|x_j-y_j|^q)+\frac{j}{q}D(|x_j-y_j|^q)\otimes
D(e^{-k(d(x_j)+d(y_j))}),$$
$$A_3:=e^{-k(d(x_j)+d(y_j))}\frac{j}{q}D^2(|x_j-y_j|^q).$$
For $A_1$ and $A_3$ we have $$A_1\leq Cj|x_j-y_j|^q \left(%
\begin{array}{cc}
  I & 0 \\
  0 & I \\
\end{array}%
\right),$$
$$ A_3\leq (q-1)j|x_j-y_j|^{q-2}\left(%
\begin{array}{cc}
  I & -I \\
  -I & I \\
\end{array}%
\right).$$ Here and henceforth, as usual, the letter $C$ denotes
various constants independent of $j$. Now we consider the quantity
$\langle A_2(\xi,\eta),(\xi,\eta)\rangle$ for $\xi,\,\eta\in\R^N$.
We have
\begin{equation*}\begin{split}\langle
A_2(\xi,\eta),(\xi,\eta)\rangle&=2kj|x_j-y_j|^{q-2}e^{-k(d(x_j)+d(y_j))}[\langle
Dd(x_j)\otimes(x_j-y_j)(\eta-\xi),\xi\rangle\\&+\langle
Dd(y_j)\otimes(x_j-y_j)(\eta-\xi),\eta\rangle]\\&\leq C
j|x_j-y_j|^{q-1}|\xi-\eta|(|\xi|+|\eta|)\\&\leq C
j|x_j-y_j|^{q-1}\left(\frac{|\xi-\eta|^2}{|x_j-y_j|}+\frac{(|\xi|+|\eta|)^2}{4}|x_j-y_j|\right)\\&\leq
C \left[j|x_j-y_j|^{q-2}|\xi-\eta|^2+
j|x_j-y_j|^{q}(|\xi|^2+|\eta|^2)\right].\end{split}\end{equation*}
The last inequality can be rewritten  equivalently in this
way$$A_2\leq C
j|x_j-y_j|^{q-2}\left(%
\begin{array}{cc}
  I & -I \\
  -I & I \\
\end{array}%
\right)+C j|x_j-y_j|^{q}\left(%
\begin{array}{cc}
  I & 0 \\
  0 & I \\
\end{array}%
\right).$$ Finally if we choose $$\epsilon
=\frac{1}{j|x_j-y_j|^{q-2}},$$ we get the same estimates for the
matrix $\epsilon(D^2\phi(x_j,y_j))^2$. In conclusion we have
\begin{equation*}\begin{split}D^2\phi(x_j,y_j)+\epsilon(D^2\phi(x_j,y_j))^2&\leq C
j|x_j-y_j|^{q-2}\left(%
\begin{array}{cc}
  I & -I \\
  -I & I \\
\end{array}%
\right)\\& +C j|x_j-y_j|^{q}\left(%
\begin{array}{cc}
  I & 0 \\
  0 & I \\
\end{array}%
\right),\end{split}\end{equation*}and \eqref{tm1lemis} implies
\eqref{claimmatrix}. The Claim 2 is proved.

Now, multiplying the inequality \eqref{claimmatrix} on the left
for the non-negative symmetric matrix
$$\left(%
\begin{array}{cc}
  \sigma(p_j)\sigma(p_j)&\sigma(p_j)\sigma(r_j) \\
\sigma(r_j)\sigma(p_j)& \sigma(r_j)\sigma(r_j) \\
\end{array}%
\right)=\left(%
\begin{array}{cc}\sigma(p_j)&\sigma(p_j)\sigma(r_j) \\
\sigma(r_j)\sigma(p_j)& \sigma(r_j) \\
\end{array}%
\right),$$  taking traces and using \eqref{infp-infq} and
\eqref{graddiff0}, we get
\begin{equation*}\begin{split}\tr
(\sigma(p_j)(X_j-\widetilde{X_j}))-\tr(\sigma(r_j)(Y_j-\widetilde{Y_j}))
&\leq \zeta_j\tr[(\sigma(p_j)-\sigma(r_j))^2]\leq
\frac{8\zeta_j}{|p_j|^2}|p_j-r_j|^2\\&\leq
C\frac{j|x_j-y_j|^{q-2}j^2|x_j-y_j|^{2q}}{j^2|x_j-y_j|^{2(q-1)}}\\&=
Cj|x_j-y_j|^q.\end{split}\end{equation*}

Now using  that $u$ and $w$ are respectively sub and supersolution
we compute
\begin{equation*}\begin{split}
-(\lam+c(x_j))u(x_j)&\leq \tr(\sigma(p_j)X_j)+b(x_j)\cdot p_j
\\&\leq \tr(\sigma(p_j)(X_j-\widetilde{X_j}))+
b(x_j)\cdot p_j + O\left(j|x_j-y_j|^{q}\right)
\\&\leq  \tr(\sigma(r_j)(Y_j-\widetilde{Y_j}))+b(x_j)\cdot
p_j+ O\left(j|x_j-y_j|^{q}\right)\\& \leq -(\tau+c(y_j))w(y_j)+
b(x_j)\cdot p_j-b(y_j)\cdot r_j+ O\left(j|x_j-y_j|^{q}\right).
\end{split}\end{equation*}
The quantity $b(x_j)\cdot p_j-b(y_j)\cdot r_j$ goes to 0 as
$j\rightarrow+\infty$. Indeed, since $m>0$ and $w$ is positive and
bounded, the estimate \eqref{stimau-v1} of Lemma \ref{regolarita}
holds for $u$ and $w$; using it in \eqref{dis1} and dividing by
$|x_j-y_j|\neq 0$ we obtain
\begin{equation*}\frac{j}{q}|x_j-y_j|^{q-1}\leq
C_0e^{k(d(x_j)+d(y_j))}\leq C.\end{equation*}Then by
\eqref{graddiff0} we conclude that the sequences $\{p_j\}$ and
$\{r_j\}$ are bounded, so that, since in addition $|p_j-r_j|\leq
Cj|x_j-y_j|^{q}\rightarrow0$ as $j\rightarrow+\infty$, up to
subsequence $p_j,r_j\rightarrow p_0$ as $j\rightarrow+\infty$.

Hence, sending $j\rightarrow+\infty$ we obtain
\begin{equation*}-(\lam+c(\zs))u(\zs)\leq
-(\tau+c(\zs))w(\zs).\end{equation*}

If $\tau+c(\zs)>0$, using \eqref{w(x)} we get
$$-(\lam+c(\zs))u(\zs)\leq-(\tau+c(\zs))\frac{\gamma}{\gamma'}u(\zs),$$
and taking $\gamma$ sufficiently close to $\gamma'$ in order that
$\frac{\tau\frac{\gamma}{\gamma'}-\lam}{1-\frac{\gamma}{\gamma'}}>
|c|_\infty,$ we obtain a contradiction. Finally if
$\tau+c(\zs)\leq0$ we have
$$-(\lam+c(\zs))u(\zs)\leq
-(\tau+c(\zs))w(\zs)\leq -(\tau+c(\zs))u(\zs),$$once more a
contradiction since $\lam<\tau$. \finedim

\dims  \textbf{of Lemma \ref{lemxdivy}.} Without loss of
generality we can assume that $\xs=0$.

Since the minimum is strict there exists a small $\delta>0$ such
that $$v(0)<v(x)+C|x|^qe^{-kd(x)}\quad \text{for any
}x\in\Oms,\,0<|x|\leq\delta.$$ Since $v$ is not locally constant
and $q>1$ for any $n>\delta^{-1}$ there
 exists $(t_n,z_n)\in B(0,\frac{1}{n})^2\cap\Oms^2$ such that
 $$v(t_n)>v(z_n)+C|z_n-t_n|^qe^{-kd(z_n)}.$$
 Consequently, for $n>\delta^{-1}$ the minimum of the function
 $v(x)+C|x-t_n|^qe^{-kd(x)}$ in $\overline{B}(0,\delta)\cap\Oms$ is not achieved on $t_n$. Indeed
 $$\min_{|x|\leq\delta,\,x\in\Oms}(v(x)+C|x-t_n|^qe^{-kd(x)})\leq
 v(z_n)+C|z_n-t_n|^qe^{-kd(z_n)}<v(t_n).$$
Let $y_n\neq t_n$ be some point in $\overline{B}(0,\delta)\cap
\Oms$ on which the minimum is achieved. Passing to the limit as
$n$ goes to infinity, $t_n$ goes to 0 and, up to subsequence,
$y_n$ converges to some $y\in \overline{B}(0,\delta)\cap\Oms$. By
the lower semicontinuity of $v$ and the fact that 0 is a local
minimum of $v(x)+C|x|^qe^{-kd(x)}$ we have
$$v(0)\leq v(y)+C|y|^qe^{-kd(y)}\leq\liminf_{n\rightarrow +\infty}(v(y_n)+C|y_n|^qe^{-kd(y_n)}),$$
and using that $v(0)+C|t_n|^qe^{-kd(0)}\geq
v(y_n)+C|y_n-t_n|^qe^{-kd(y_n)},$ one has
$$v(0)\geq\limsup_{n\rightarrow +\infty}(v(y_n)
+C|y_n|^qe^{-kd(y_n)}).$$ Then
$$v(0)=v(y)+C|y|^qe^{-kd(y)}=\lim_{n\rightarrow +\infty}(v(y_n)+C|y_n|^qe^{-kd(y_n)}).$$
Since 0 is a strict local minimum of $v(x)+C|x|^qe^{-kd(x)}$, the
last equalities imply that $y=0$ and $v(y_n)$ goes to $v(0)$ as
$n\rightarrow +\infty$. Then for large $n$, $y_n$ is an interior
point of $B(0,\delta)$ so that the function
$$\varphi(x)=v(y_n)+C|y_n-t_n|^qe^{-kd(y_n)}-C|x-t_n|^qe^{-kd(x)}$$ is a test function for $v$ at $y_n$. Moreover, the gradient of $\varphi$
$$D\varphi(x)=-Cq|x-t_n|^{q-2}e^{-kd(x)}(x-t_n)+kC|x-t_n|^qe^{-kd(x)}Dd(x)$$
is different from 0 at $x=y_n$ for small $\delta$, indeed
$$|D\varphi(y_n)|\geq
C|y_n-t_n|^{q-1}e^{-kd(y_n)}(q-k|y_n-t_n|)\geq
C|y_n-t_n|^{q-1}e^{-kd(y_n)}(q-2k\delta)>0.$$

 Using
\eqref{sferaest}, if $y_n\in\partial \Om$ we have
\begin{equation*}\begin{split}\langle
D\varphi(y_n),\overrightarrow{n}(y_n)\rangle\leq
C|y_n-t_n|^q\left(\frac{q}{2r}-k\right)<0,
\end{split}\end{equation*} since $k>q/(2r)$. Then we conclude that
$$\tr\left(\sigma(D\varphi(y_n))D^2\varphi(y_n)\right)+b(y_n)\cdot
D\varphi(y_n)-\beta(v(y_n))\leq g(y_n).$$ Observe that
$D^2\varphi(y_n)=|y_n-t_n|^{q-2}M,$ where $M$ is a bounded matrix.
Hence, from the last inequality we get
\begin{equation*}C_0|y_n-t_n|^{q-2}-\beta(v(y_n))\leq
g(y_n),\end{equation*}for some constant $C_0$. Passing to the
limit, since $\beta$ and $g$ are upper semicontinuous we obtain
$$-\beta (v(0))\leq g(0),$$ which is the desired conclusion.
\finedim We conclude  sketching the proof of Theorem
\ref{dirichletcomp}.

 \dims \textbf{of Theorem
\ref{dirichletcomp}.} Suppose by contradiction that
$\max_{\Oms}(u-v)=m>0$. Since $u\leq v$ on the boundary, the
supremum is achieved inside $\Om$. Let us define for $j\in\N$ and
some $q>2$
$$\psi(x,y)=u(x)-v(y)-\frac{j}{q}|x-y|^q.$$ Suppose that
$(x_j,y_j)$ is a maximum point for $\psi$ in $\Oms^2$. Then
$|x_j-y_j|\rightarrow0$ as $j\rightarrow+\infty$ and up to
subsequence $x_j,y_j\rightarrow\xs$, $u(x_j)\rightarrow u(\xs)$,
$v(y_j)\rightarrow v(\xs)$ and $j|x_j-y_j|^q\rightarrow0$ as
$j\rightarrow+\infty$. Moreover, $\xs$ is such that
$u(\xs)-v(\xs)=m$ and we can choose $x_j\neq y_j$. Recalling by
Remark \ref{stimadirich} that the estimate \eqref{stimau-v1} holds
in $\Om$, we can proceed as in the proof of Theorem \ref{maxpneum}
to get $$-c(\xs) u(\xs)\leq -c(\xs) v(\xs).$$ This is a
contradiction since $c(\xs)<0$. \finedim

 \subsection{The Maximum Principle for $c(x)$ chan\-ging sign: an example.}
In the previous subsections we have proved that $\di+ b(x) \cdot
D+ c(x)$ with the Neumann boundary condition satisfies the maximum
principle if $c(x)\leq 0$ or without condition on the sign of
$c(x)$ provided $\lams>0$. In this subsection we want to show that
these two cases do not coincide, i.e., that there exists some
$c(x)$ which changes sign in $\Om$ such that the associated
principal eigenvalue $\lams$ is positive. To prove this, by
definition of $\lams$, it suffices to find a function $c(x)$
changing sign for which there exists a bounded positive
supersolution of
\begin{equation}\label{equazsottosol}
\begin{cases}
 \di v+b(x)\cdot Dv+(c(x)+\lam)v = 0 & \text{in} \quad\Om \\
 \dfrac{\partial v}{\partial \overrightarrow{n}}= 0 & \text{on} \quad\partial\Om, \\
 \end{cases}
 \end{equation}for some $\lam>0$. For simplicity, let us suppose that
$b\equiv 0$ and $\Om$ is the ball of center 0 and radius R. We
will look for $c$ such that:

\begin{equation}\label{cprop}
\begin{cases}
c(x)<0 &\text{if } R-\epsilon<|x|\leq R\\
c(x)\leq -\beta_1 & \text{if }\rho<|x|\leq R-\epsilon\\
c(x)\leq \beta_2 & \text{if }|x|\leq \rho,\\
 \end{cases}
 \end{equation}
 where $0<\rho<R-\epsilon$ and $\epsilon$, $\beta_1,\, \beta_2$ are positive constants. Remark that in the ball of radius $\rho$,
$c(x)$ may assume positive values. Following \cite{p1}, it is
possible to construct a supersolution of \eqref{equazsottosol} if
$\epsilon$ is small enough and
$$\beta_2<\frac{k^2e^{-k\rho}}{\frac{k}{4}(R-\rho)+\frac{2k}{\beta_1(R-\rho)}+1-e^{-k\rho}},$$
for some $k>0$. From the last relation we can see that choosing
$k=\frac{1}{\rho}$ the term on the right-hand side goes to
$+\infty$ as $\rho\rightarrow 0^+$, that is, if the set where
$c_0(x)$ is positive becomes smaller then the values of $c_0(x)$
in this set can be very large. On the contrary, for any value of
$k$, if $\rho\rightarrow R^-$ then $\beta_2$ goes to 0. Finally
for any $k$ if $\beta_1\rightarrow 0^+$, then again $\beta_2$ goes
to 0.
\section{Some existence results} This section is devoted to the
problem of the existence of a solution of
\begin{equation}\label{existence}
\begin{cases}
 \di u+ b(x)\cdot Du+(c(x)+\lam) u=  g(x) & \text{in} \quad\Om \\
 \dfrac{\partial u}{\partial \overrightarrow{n}}= 0 & \text{on} \quad\partial\Om.\\
 \end{cases}
 \end{equation}
 The first existence result for \eqref{existence} is obtained when
 $\lam=0$ and $c<0$, via Perron's method. Then, we will prove the
 existence of a positive solution of \eqref{existence} when $g$
 is  non-positive and $\lam<\lams$ (without condition on the sign of $c$).
  These two results will allow us to prove that the Neumann problem  \eqref{existence}
 is solvable for any right-hand side if $\lam<\lams$. Finally,  we will prove the existence of a positive principal
 eigenfunction corresponding to $\lams$, that is a solution of
 \eqref{existence} when $g\equiv0$ and $\lam=\lams$.

 Comparison results guarantee for \eqref{existence} the uniqueness of the solution
 when $c<0$ and when $\lam<\lams$ and $g<0$ or $g>0$.

\begin{thm}\label{princonfc<0}Suppose that $c<0$ and $g$ is continuous
on $\Oms$. If $u\in USC(\Oms)$ and $v\in LSC(\Oms)$ are
respectively  viscosity sub and supersolution of
\begin{equation}\label{princonfl<ls2}
\begin{cases}
 \di u+b(x)\cdot Du+c(x) u=  g(x) & \text{in} \quad\Om \\
 \dfrac{\partial u}{\partial \overrightarrow{n}}= 0 & \text{on} \quad\partial\Om, \\
 \end{cases}
 \end{equation}
 with $u$ and $v$ bounded or $v\geq0$ and bounded, then $u\leq v$ on
$\Oms.$ Moreover \eqref{princonfl<ls2} has a unique viscosity
solution.
\end{thm}
\dim  We suppose by contradiction that $\max_{\Oms}(u-v)=m>0$.
Repeating the proof of Theorem \ref{maxpneum} taking $v$ as $w$,
we arrive to the following inequality
$$-c(\zs) u(\zs)\leq -c(\zs)
v(\zs),$$ where $\zs\in\Oms$ is such that $u(\zs)-v(\zs)=m>0$.
This is a contradiction since $c(\zs)<0$.

The existence of a solution follows from Perron's method of Ishii,
see e.g. \cite{cil}, and the comparison result just proved,
provided there is a bounded subsolution and a bounded
supersolution of \eqref{princonfl<ls2}. Since $c$ is negative and
continuous on $\Oms$, there exists $c_0>0$ such that $c(x)\leq
-c_0$ for every $x\in\Oms$. Then
$$u_1:=-\frac{|g|_\infty}{c_0},\quad u_2:=\frac{|g|_\infty}{c_0}$$ are
respectively a bounded sub and supersolution of
\eqref{princonfl<ls2}.

Define
$$u(x):=\sup\{\varphi(x)|\,u_1\leq
\varphi\leq u_2\text{ and }\varphi \text{ is a subsolution of
\eqref{princonfl<ls2} } \},$$we claim that $u$ is a solution of
\eqref{princonfl<ls2}. We first show that the upper semicontinuous
envelope of u defined as $$u^*(x):=\lim_{ \rho\downarrow
0}\,\sup\{u(y):y\in\Oms\text{ and } |y-x|\leq \rho\}$$ is a
subsolution of \eqref{princonfl<ls2}. Indeed if $(p,X)\in
J^{2,+}u(x_0)$ and $p\neq 0$ then by the standard arguments of the
Perron's method it can be proved that tr$(\sigma(p)X)+b(x_0)\cdot
p+c(x_0)u(x_0)\geq g(x_0)$ if $x_0\in\Om$ and
$(-\tr(\sigma(p)X)-b(x_0)\cdot p-c(x_0)u(x_0)+g(x_0))\wedge
\langle p,\overrightarrow{n}(x_0)\rangle\leq 0$ if
$x_0\in\partial\Om$.

Now suppose $u^*\equiv k$ in a neighborhood of $x_0\in\Oms$. If
$x_0\in\partial \Om$ clearly $u^*$ is subsolution at $x_0$. Assume
that $x_0$  is an interior point of $\Om$. We may choose a
sequence of subsolutions $(\varphi_n)_n$ and a sequence of points
$(x_n)_n$ in $\Om$ such that $x_n\rightarrow x_0$ and
$\varphi_n(x_n)\rightarrow k$. Suppose that $|x_n-x_0|<a_n$ with
$a_n$ decreasing to 0 as $n\rightarrow+\infty$. If, up to
subsequence, $\varphi_n$ is constant in $B(x_0,a_n)$ for any $n$,
then passing to the limit in the relation $c(x_n)
\varphi_n(x_n)\geq g(x_n)$ we get $c(x_0) k\geq g(x_0)$ as
desired. Otherwise, suppose that for any $n$ $\varphi_n$ is not
constant in $B(x_0,a_n)$. Repeating the argument of Lemma
\ref{lemxdivy} we find a sequence
$\{(t_n,y_n)\}_{n\in\N}\subset\Om^2$ and a small $\delta>0$ such
that $|t_n-x_0|<a_n$, $|y_n-x_0|\leq \delta$, $t_n\neq y_n$,
$\varphi_n(x)-|x-t_n|^q\leq\varphi_n(y_n)-|y_n-t_n|^q$ for any
$x\in B(x_0,\delta)$, with $q>2$ and $u^*\equiv k$ in
$\overline{B}(x_0,\delta)$. Up to subsequence $y_n\rightarrow y\in
\overline{B}(x_0,\delta)$ as $n\rightarrow+\infty$. We have
\begin{equation*}\begin{split}k& =\lim_{n\rightarrow+\infty}(\varphi_n(x_n)-|x_n-t_n|^q)\leq
\liminf_{n\rightarrow+\infty}(\varphi_n(y_n)-|y_n-t_n|^q)\\&\leq\limsup_{n\rightarrow+\infty}(\varphi_n(y_n)-|y_n-t_n|^q)\leq
k-|y-x_0|^q.\end{split}\end{equation*} The last inequalities imply
that $y=x_0$ and $\varphi_n(y_n)\rightarrow k$. Then, for large
$n$, $y_n$ is an interior point of $B(x_0,\delta)$ and
$\phi_n(x):= \varphi_n(y_n)-|y_n-t_n|^q + |x-t_n|^q$ is a test
function for $\varphi_n$ at $y_n$. Passing to the limit as
$n\rightarrow+\infty$ in the relation $\di\phi_n(y_n)+b(y_n)\cdot
D\phi_n(y_n)+c(y_n)\varphi_n(y_n)) \geq g(y_n)$, we get again
$c(x_0) k\geq g(x_0)$. In conclusion $u^*$ is a subsolution of
\eqref{princonfl<ls2}. Since $u^*\leq u_2$, it follows from the
definition of $u$ that $u=u^*$.

Finally the lower semicontinuous envelope of $u$ defined as
$$u_*(x):=\lim_{ \rho\downarrow0}\,\inf\{u(y):y\in\Oms\text{ and }
|y-x|\leq \rho\}$$is a supersolution. Indeed, if it is not, the
Perron's method provides  a viscosity subsolution of
\eqref{princonfl<ls2} greater than $u$, contradicting the
definition of $u$. If $u_*\equiv k$ in a neighborhood of
$x_0\in\Om$ and $c(x_0) k>g(x_0)$ then for small $\delta$ and
$\rho$, the subsolution is
\begin{equation*}u_{\delta,\rho}(x):=
\begin{cases}\max\{u(x),k+\frac{\delta \rho^2}{8}-\delta|x-x_0|^2\} &\text{if }|x-x_0|<\rho,\\
u(x)&\text{otherwise}.\end{cases}
\end{equation*} Hence
$u_*$ is a supersolution of \eqref{princonfl<ls2} and then, by
comparison, $u^*=u\leq u_*$, showing that $u$ is continuous and is
a solution.

The uniqueness of the solution is an immediate consequence of the
comparison principle just proved. \finedim
\begin{thm}\label{comparisonl<ls}
Suppose $g\in LSC(\Oms)$, $h\in USC(\Oms)$, $h\leq0$, $h\leq g$
and $g(x)>0$ if $h(x)=0$. Let $u\in USC(\Oms)$ be a viscosity
subsolution of \eqref{existence} and $v\in LSC(\Oms)$ be a bounded
positive viscosity supersolution of \eqref{existence} with $g$
replaced by $h$. Then $u\leq v$ on $\Oms.$
 \begin{rem}{\em The existence of a such $v$ implies
 $\lam\leq\lams$.}\end{rem}
 \end{thm}
 \dim
 It suffices to prove the theorem for $h<g$. Indeed, for $l>1$ the function $lv$ is a
 supersolution of \eqref{existence} with right-hand side $lh(x)$ and by the assumptions on $h$ and $g$, $lh<g$. If $u\leq lv$
for any $l>1$, passing to the limit as $l\rightarrow1^+$, one
obtains $u\leq v$ as desired.

Hence we can assume $h<g$. By upper semicontinuity
$\max_{\Oms}(h-g)=-M<0$. Suppose by contradiction that $u>v$
somewhere in $\Oms$. Then there exists $\ys\in\Oms$ such that
$$\gamma':=\frac{u(\ys)}{v(\ys)}=\max_{x\in\Oms}\frac{u(x)}{v(x)}>1.$$ Define
$w=\gamma v$ for some $1\leq\gamma<\gamma'$. Since $h\leq 0$ and
$\gamma\geq1$, $\gamma h\leq h$ and then $w$ is still a
supersolution of \eqref{existence} with right-hand side $h$. The
supremum of $u-w$ is strictly positive then, by upper
semicontinuity, there exists $\xs\in\Oms$ such that
$u(\xs)-w(\xs)=\max_{\Oms}(u-w)>0$. We have $u(\xs)>w(\xs)$ and
$w(\xs)\geq\frac{\gamma}{\gamma'} u(\xs)$. Repeating the proof of
Theorem \ref{maxpneum}, we get
$$g(\zs)-(\lam+c(\zs))u(\zs)\leq
h(\zs)-(\lam+c(\zs))w(\zs),$$where $\zs$ is some point in $\Oms$
where the maximum of $u-w$ is attained. If $\lam+c(\zs)\leq 0$,
then
$$-(\lam+c(\zs))u(\zs)\leq h(\zs)-g(\zs)
-(\lam+c(\zs))w(\zs)<-(\lam+c(\zs))u(\zs),$$ which is a
contradiction. If $\lam+c(\zs)>0$, then
$$-(\lam+c(\zs))u(\zs)\leq
h(\zs)-g(\zs)-(\lam+c(\zs))\frac{\gamma}{\gamma'}u(\zs).$$ If we
choose $\gamma$ sufficiently close to $\gamma'$ in order that
$$|\lam+c|_\infty\left(\frac{\gamma}{\gamma'}-1\right)\max_{\Oms}u\geq -\frac{M}{2},$$ we get
 once more a contradiction. \finedim
\begin{thm}\label{esistl<ls}
Suppose that $\lam<\lams$, $g\leq 0$, $g\not\equiv0$ and $g$ is
continuous on $\Oms$, then there exists a positive viscosity
solution of \eqref{existence}. If $g<0$, the solution is unique.
\end{thm}
\dim We follow the proof of Theorem 7 of \cite{bd}.

If $\lam<-|c|_\infty$ then the existence of the solution is
guaranteed by Theorem \ref{princonfc<0}. Let us suppose $\lam\geq
-|c|_\infty$ and define by induction the sequence $(u_n)_n$ by
$u_1=0$ and $u_{n+1}$ as the solution of
\begin{equation*}
\begin{cases}
 \di u_{n+1}+b(x)\cdot Du_{n+1}+(c(x)-|c|_\infty-1)u_{n+1}
  =g-(\lam+|c|_\infty+1)u_n & \text{in} \quad\Om \\
 \dfrac{\partial u_{n+1}}{\partial \overrightarrow{n}}= 0 & \text{on  } \partial\Om, \\
 \end{cases}
 \end{equation*}
which exists by Theorem \ref{princonfc<0}. By the comparison
principle, since $g\leq 0$ and $g\not\equiv 0$ the sequence is
positive and increasing.

We claim that $(u_n)_n$ is also bounded. Suppose that it is not,
then dividing by $|u_{n+1}|_\infty$ and defining
$v_n:=\frac{u_n}{|u_n|_\infty}$ one gets that $v_{n+1}$ is a
solution of
\begin{equation*}
\begin{cases}
 \di v_{n+1}+b(x)\cdot Dv_{n+1}+(c(x)-|c|_\infty-1)v_{n+1}
 \\ \quad=\frac{g}{|u_{n+1}|_\infty}    -(\lam+|c|_\infty+1)\frac{u_n}{|u_{n+1}|_\infty}& \text{in} \quad\Om \\
 \dfrac{\partial v_{n+1}}{\partial \overrightarrow{n}}= 0 & \text{on} \quad\partial\Om. \\
 \end{cases}
 \end{equation*}
 By Corollary \ref{corcomp},  $(v_n)_n$ converges to a positive function $v$ with $|v|_\infty=1$, which
 satisfies
 \begin{equation*}
\begin{cases}
 \di v+b(x)\cdot Dv+(c(x)+\lam) v\\\quad=(\lam+|c|_\infty+1)(1-k)v\geq 0& \text{in} \quad\Om \\
 \dfrac{\partial v_{n+1}}{\partial \overrightarrow{n}}= 0 & \text{on} \quad\partial\Om, \\
 \end{cases}
 \end{equation*}where
 $k:=\lim_{n\rightarrow+\infty}\frac{|u_n|_\infty}{|u_{n+1}|_\infty}\leq
 1$. This contradicts the maximum principle, Theo\-rem
 \ref{maxpneum}.

 Then $(u_n)_n$ is bounded and  letting $n$ go to infinity, by the
compactness result, the sequence converges to a function $u$ which
is a solution. Moreover, the solution is positive on $\Oms$ by the
strong minimum principle, Proposition \ref{p1}.

If $g<0$, the uniqueness of the solution follows from Theorem
\ref{comparisonl<ls}.\finedim
\begin{rem}\label{thmslsot}{\em Clearly, since the operator $\di $ is odd, by Theorem \ref{esistl<ls}, there exists a negative solution of
 \eqref{existence} for
$\lam<\underline{\lam}$ and $g\geq 0$, $g\not\equiv0$, which is
unique if $g>0$.}
\end{rem}

\begin{thm}Suppose that
$\lam<\lams$ and  $g$ is continuous on $\Oms$, then there exists a
viscosity solution of \eqref{existence}.
\end{thm}
 \dim If $g\equiv 0$, by the maximum principle the only solution is $u\equiv 0$. Let us suppose $g\not\equiv0$.
 Since $\lam<\lams$ by Theorem
 \ref{esistl<ls} there exist
 $v_0$ positive viscosity solution of
 \eqref{existence} with right-hand side $-|g|_\infty$ and $u_0$ negative viscosity solution
 of \eqref{existence} with right-hand side $|g|_\infty$.

 Let us suppose $\lam+|c|_\infty\geq 0$. Let $(u_n)_n$ be the sequence defined in the proof of Theorem \ref{esistl<ls} with
 $u_1=u_0$, then by comparison Theorem \ref{princonfc<0} we have $u_0=u_1\leq
u_2\leq
 ...\leq v_0$. Hence, by the compactness Corollary \ref{corcomp}
 the sequence converges to a continuous function which is the
 desired solution.\finedim

\begin{thm}[Existence of principal eigenfunctions]\label{esistautof}There exists $\phi>0$ on $\Oms$
viscosity solution of
\begin{equation*}
\begin{cases}
 \di \phi+b(x)\cdot D\phi+(c(x)+\lams)\phi= 0 & \text{in} \quad\Om \\
 \dfrac{\partial \phi}{\partial \overrightarrow{n}}= 0 & \text{on} \quad\partial\Om. \\
 \end{cases}
 \end{equation*}Moreover $\phi$ is Lipschitz continuous on $\Oms$.
\end{thm}
\dim Let $\lam_n$ be an increasing sequence which converges to
$\lams$. Let $u_n$ be the positive solution of \eqref{existence}
with $\lam=\lam_n$ and $g\equiv -1$. By  Theorem \ref{esistl<ls}
the sequence $(u_n)_n$ is well
 defined. Following the argument of the proof of Theorem 8 of
 \cite{bd}, it can proved that it is unbounded, otherwise one would
 contradict the definition of $\lams$.
 Then, up to subsequence $|u_n|_\infty\rightarrow+\infty$ as
 $n\rightarrow+\infty$ and defining $v_n:=\frac{u_n}{|u_n|_\infty}$
 one gets that $v_n$ satisfies \eqref{existence} with
 $\lam=\lam_n$ and $g\equiv  -\frac{1}{|u_n|_\infty}$.
 Then by Corollary \ref{corcomp}, we can extract a subsequence converging to a positive
 function $\phi$ with $|\phi|_\infty=1$ which is the desired solution.
 By  Theorem \ref{regolaritathm} the solution is also Lipschitz continuous on
 $\Oms$. \finedim

\section{A decay estimate for solutions of the evolution problem}
In this section we want to study the asymptotic behavior as
$t\rightarrow+\infty$ of the solution $h(t,x)$ of the evolution
problem
\begin{equation}\label{sisevolu}
\begin{cases}
h_t=\di h+c(x)h & \text{in } (0,+\infty)\times\Om\\
\dfrac{\partial h}{\partial \overrightarrow{n}}=0 &\text{on } [0,+\infty)\times \partial\Om\\
h(0,x)=h_0(x) & \text{for }x\in\Om,\\
\end{cases}
 \end{equation}where $h_0$ is a continuous function on $\Oms$.
 As in \cite{j} and in \cite{jk} we use the se\-mi\-con\-ti\-nuous
 extensions of the function $(p,X)\rightarrow\tr (\sigma(p)X)$ to
 define the viscosity solutions of \eqref{sisevolu}. For $X\in\emph{S(N)}$, let us denote its smaller and larger
 eigenvalue respectively by $m(X)$ and $M(X)$, that is
 $$m(X):=\min_{|\xi|=1}\langle X\xi,\xi\rangle,$$
 $$M(X):=\max_{|\xi|=1}\langle X\xi,\xi\rangle.$$

 \begin{de} Any function $u\in USC([0,+\infty)\times\Oms)$ (resp., $u\in LSC([0,+\infty)\times\Oms)$) is
 called \emph{viscosity subsolution} (resp., \emph{supersolution}) of \eqref{sisevolu}
  if for any $x\in\Oms$, $u(0,x)\leq h_0(x)$ (resp., $u(0,x)\geq h_0(x)$) and if the following conditions hold
\begin{itemize}
\item[(i)] For every $(t_0,x_0)\in  (0,+\infty)\times\Om$, for all $\varphi\in C^2([0,+\infty)\times\Oms)$,
 such that $u-\varphi$ has a local maximum (resp., minimum)
at $(t_0,x_0)$, one has
\begin{equation*}
\begin{cases}\varphi_t(t_0,x_0)\leq\di
\varphi(t_0,x_0)+c(x_0)u(t_0,x_0)\,(\text{resp., }\geq) &\text{ if
}D\varphi(t_0,x_0)\neq 0,\\ \varphi_t(t_0,x_0)\leq
M(D^2\varphi(t_0,x_0))+c(x_0)u(t_0,x_0) &\text{ if
}D\varphi(t_0,x_0)= 0\\(\text{resp., }\varphi_t(t_0,x_0)\geq
m(D^2\varphi(t_0,x_0))+c(x_0)u(t_0,x_0)).
\end{cases}
 \end{equation*}
\item[(ii)]For every $(t_0,x_0)\in (0,+\infty)\times\partial\Om$, for all $\varphi\in C^2([0,+\infty)\times\Oms)$, such that
$u-\varphi$ has a local maximum (resp., minimum) at $(t_0,x_0)$
and $D\varphi(t_0,x_0)\neq 0$, one has
\begin{equation*}
(\varphi_t(t_0,x_0)-\di
\varphi(t_0,x_0)-c(x_0)u(t_0,x_0))\wedge\langle
D\varphi(t_0,x_0),\overrightarrow{n}(x_0)\rangle\leq 0.
 \end{equation*}(resp., $$(\varphi_t(t_0,x_0)-\di
\varphi(t_0,x_0)-c(x_0)u(t_0,x_0))\vee\langle
D\varphi(t_0,x_0),\overrightarrow{n}(x_0)\rangle\geq 0.)$$

\end{itemize}
\end{de}
Remark that if $(t_0,x_0)\in (0,+\infty)\times\partial\Om$ and
$D\varphi(t_0,x_0)= 0$, then the boundary condition is satisfied.

 We will show that if the principal eigenvalue of the stationary operator as\-so\-cia\-ted to
 \eqref{sisevolu} is positive, then $h$ decays to zero exponentially and that the
 rate of the decay depends on it. Let $\lams$ and $v$ be respectively the principal eigenvalue and a principal eigenfunction, i.e., $v$ is
 a positive
 solution of
 \begin{equation*}
\begin{cases}
\di v+(c(x)+\lams)v=0& \text{in }\Om\\
\dfrac{\partial v}{\partial \overrightarrow{n}}=0  &\text{on
}\partial\Om.
\end{cases}
 \end{equation*}

\begin{prop}Let $h\in C(\Oms\times[0,+\infty))$ be a solution of \eqref{sisevolu} then
\begin{equation}\label{crescita}\sup_{\Om\times[0,+\infty)}\frac{h(t,x)e^{\lams t}}{v(x)}\leq
\sup_{\Om}\frac{h_0^+(x)}{v(x)},\end{equation} where
$h_0^+=\max\{h_0,0\}$ denotes the positive part of $h_0$.
\end{prop}\dim It suffices to prove that, fixed $\lam<\lams$
$$\sup_{[0,T)\times\Om}\frac{h(t,x)e^{\lam t}}{v(x)}\leq
\sup_{\Om}\frac{h_0^+(x)}{v(x)},$$ for any $T>0$. This implies
that
$$\sup_{[0,T)\times\Om}\frac{h(t,x)e^{\lams t}}{v(x)}\leq
\sup_{\Om}\frac{h_0^+(x)}{v(x)},$$  for any $T>0$ and consequently
\eqref{crescita}. Let us denote $H(t,x)=h(t,x)e^{\lam t}$, it is
easy to see that $H(t,x)$ satisfies
\begin{equation}\label{sish}
\begin{cases}
H_t=\di H+(c(x)+\lam)H & \text{in } [0,+\infty)\times\Om\\
\dfrac{\partial H}{\partial \overrightarrow{n}}=0 &\text{on } [0,+\infty)\times\partial\Om\\
H(0,x)=h_0(x) & \text{for }x\in\Om.\\
\end{cases}
 \end{equation}Suppose by contradiction that there exists $T>0$
 such that \begin{equation}\label{contrevol}\gamma':=\sup_{[0,T)\times \Om}\frac{h(t,x)e^{\lam
t}}{v(x)}> \sup_{\Om}\frac{h_0^+(x)}{v(x)}=:\overline{h}\geq
0.\end{equation} Let us denote $w=\gamma v$, where
$$\overline{h}<\gamma<\gamma'$$ and $\gamma$ is sufficiently close
to $\gamma'$ in order that
\begin{equation}\label{sceltagamma}\frac{\lams\frac{\gamma}{\gamma'}-\lam}{1-\frac{\gamma}{\gamma'}}>
|c|_\infty.\end{equation}Since $\gamma<\gamma'$, the function
$H-w$ has a positive maximum on $[0,T]\times\Oms$.

 Fix $q>2$, $k>\frac{q}{2r}$ and $\ep>0$ small, for $j\in\N$ we
 define the function
$$\phi(t,x,s,y)=\left(\frac{j}{q}|x-y|^q+\frac{j}{2}|t-s|^2\right)e^{-k(d(x)+d(y))}+\frac{\ep}{T-t},$$
and we consider the supremum of
$$H(t,x)-w(y)-\phi(t,x,s,y)$$ over
$([0,T)\times\Oms)^2$. Let $(t_j,x_j,s_j,y_j)$ be a point in
$(\Oms\times [0,T))^2$ where the maximum is attained. From
$$H(t_j,x_j)-w(y_j)-\phi(t_j,x_j,t_j,y_j)\leq H(t_j,x_j)-w(y_j)-\phi(t_j,x_j,s_j,y_j)$$ we deduce that
$$t_j=s_j.$$ Let $(\widehat{t},\widehat{x})\in [0,T[\times\Oms$ be
such that $H(\widehat{t},\widehat{x})-w(\widehat{x})=l>0,$ then
 for $\ep$ small enough we have
$$\frac{l}{2}\leq
H(\widehat{t},\widehat{x})-w(\widehat{x})-\frac{\ep}{T-\widehat{t}}\leq
H(t_j,x_j)-w(y_j)-\frac{\ep}{T-t_j}-
\frac{j}{q}|x_j-y_j|^qe^{-k(d(x_j)+d(y_j))}.$$ Since
$\frac{\ep}{T-t}\rightarrow+\infty$ as $t\uparrow T$, the previous
inequality implies that, up to subsequence
$(t_j,x_j,y_j)\rightarrow (\overline{t},\xs,\xs)$ as
$j\rightarrow+\infty$ with $\overline{t}<T$ and that
\begin{equation}\label{contbasso}
H(\overline{t},\xs)-w(\xs)>0.\end{equation}Moreover
$$\lim_{j\rightarrow+\infty}\frac{j}{q}|x_j-y_j|^q=0,$$ and from
\eqref{contrevol} we deduce that
\begin{equation}\label{contbasso2}w(\xs)\geq
\frac{\gamma}{\gamma'}H(\overline{t},\xs).\end{equation}Finally,
since $\gamma>\overline{h}$, it is $\overline{t}>0$. Hence for $j$
large enough, $0<t_j<T$.

As in Theorem \ref{maxpneum} the following holds true.

\textbf{Claim }\emph{ For $j$ large enough, we can choose $x_j\neq
y_j$}.

Indeed, suppose that $x_j=y_j$, then $(t_j,x_j)$ is a maximum
point for
$$U(t,x):=H(t,x)-\frac{\ep}{T-t}-e^{-kd(x_j)}\left(\frac{j}{q}|x-x_j|^q+\frac{j}{2}|t-t_j|^2\right)e^{-kd(x)},$$
and a minimum point for
$$W(t,x):=w(x)+e^{-kd(x_j)}\left(\frac{j}{q}|x-x_j|^q+\frac{j}{2}|t-t_j|^2\right)e^{-kd(x)}.$$

We prove that $(t_j,x_j)$ is not both a strict local maximum and a
strict local minimum. Indeed, in that case, if
$H(t,x)-\frac{\ep}{T-t}$ is not locally constant around
$(t_j,x_j)$, following the proof of Lemma \ref{lemxdivy}, we can
construct sequences $(t_n,x_n)_n$, $(s_n,y_n)_n$ converging to
$(t_j,x_j)$ as $n\rightarrow+\infty$, such that $(t_n,x_n)\neq
(s_n,y_n)$ and
\begin{equation*}\begin{split}\varphi(t,x)&:=C\left(\frac{|x-x_n|^q}{q}+\frac{|t-t_n|^2}{2}\right)e^{-kd(x)}+\frac{\ep}{T-t}+H(s_n,y_n)\\
&-\frac{\ep}{T-s_n}-C\left(\frac{|y_n-x_n|^q}{q}+\frac{|s_n-t_n|^2}{2}\right)e^{-kd(y_n)}\end{split}\end{equation*}is
a test function for $H(t,x)$ at $(s_n,y_n)$, where
$C=je^{-kd(x_j)}.$ If $y_n\in\partial\Om$, then $$\langle
D\varphi(s_n,y_n),\overrightarrow{n}(y_n)\rangle\geq
C\left[\left(\frac{k}{q}-\frac{1}{2r}\right)|x_n-y_n|^q+\frac{k}{2}|s_n-t_n|^2\right]>0.$$Then
$D\varphi(s_n,y_n)\neq 0$ and by definition of subsolution
$$\frac{\ep}{(T-s_n)^2}+Ce^{-kd(y_n)}(s_n-t_n)\leq \di(\varphi(s_n,y_n))+(c(y_n)+\lam)H(s_n,y_n).$$
If $y_n$ is an interior point and $D\varphi(s_n,y_n)\neq 0$, then
again the previous inequality holds true, otherwise if
 $D\varphi(s_n,y_n)= 0$, we have
$$\frac{\ep}{(T-s_n)^2}+Ce^{-kd(y_n)}(s_n-t_n)\leq
M(D^2\varphi(s_n,y_n))+(c(y_n)+\lam)H(s_n,y_n).$$ Passing to the
limit as $n\rightarrow+\infty$, from both the previous relations
we get
$$\frac{\ep}{(T-t_j)^2}\leq (c(x_j)+\lam)H(t_j,x_j).$$By definition of subsolution, we get the same inequality if $H(t,x)-\frac{\ep}{T-t}$
is locally constant around $(t_j,x_j).$

Proceeding in the same way, if either $w$ is locally constant
around $x_j$ or not, since $(t_j,x_j)$ is a strict local minimum
of $W(t,x)$, we get
$$(c(x_j)+\lams)w(x_j)\leq 0.$$ Then, passing to the limit as
$j\rightarrow+\infty$, we finally obtain
\begin{equation}\label{ultima}(c(\xs)+\lams)w(\xs)<\frac{\ep}{(T-\overline{t})^2}\leq
(c(\xs)+\lam)H(\overline{t},\xs),\end{equation} which contradicts
\eqref{sceltagamma}, \eqref{contbasso} and \eqref{contbasso2}.

Hence $(t_j,x_j)$ cannot be both a strict local maximum and a
strict local minimum. In the first case, there exists
$(s_j,y_j)\neq (t_j,x_j)$ such that
\begin{equation*}\begin{split}H(s_j,y_j)-w(x_j)-\frac{\ep}{T-s_j}-\left(\frac{j}{q}|x_j-y_j|^q+\frac{j}{2}|t_j-s_j|^2\right)e^{-k(d(x_j)+d(y_j))}
\\=H(t_j,x_j)-w(x_j)-\frac{\ep}{T-t_j}=\sup_{([0,T)\times\Om)^2}(H(t,x)-w(y)-\phi(t,x,s,y)).\end{split}\end{equation*}
As before we get that $s_j=t_j$, then $x_j\neq y_j$ and this
concludes the claim.

From the claim we deduce that $D_x\phi(t_j,x_j,t_j,y_j)$ and
$D_y\phi(t_j,x_j,t_j,y_j)$ are different from 0. Moreover there
exist $X_j,Y_j\in\emph{S(N)}$ satisfying \eqref{claimmatrix} such
that
$\left(\frac{\ep}{(T-t_j)^2},D_x\phi(t_j,x_j,t_j,y_j),X_j\right)\in\mathcal{P}^{2,+}H(t_j,x_j)$
and $(-D_y\phi(t_j,x_j,t_j,y_j),Y_j)\in J^{2,-}w(y_j)$. Now we can
proceed as in the proof of Theorem \ref{maxpneum} to obtain
\eqref{ultima} and hence to reach a contradiction.

\finedim

%$$t_j-s_j+\frac{\ep}{(T-t_j)^2}\leq
% Moreover by Theorem 3.2 of \cite{cil}
%there exist $X_j,Y_j\in\emph{S(N)}$ such that
%$$(j(t_j-s_j),D_x\phi(x_j,y_j), X_j)\in
%\mathcal{\overline{P}}^{2,+}H(x_j,t_j),$$
%$$(j(t_j-s_j),-D_y\phi(x_j,y_j), Y_j)\in
%\mathcal{\overline{P}}^{2,-}w(y_j,s_j),$$
% \begin{equation*}-\left(\frac{1}{\epsilon}+\|D^2\phi(x_j,y_j)\|\right)\left(%
%\begin{array}{cc}
%  I & 0 \\
%  0 & I \\
%\end{array}%
%\right)\leq \left(%
%\begin{array}{cc}
%  X_j & 0 \\
%  0 & -Y_j \\
%\end{array}%
%\right)\leq D^2\phi(x_j,y_j)+\epsilon (D^2\phi(x_j,y_j))^2.
%\end{equation*}
%Using the facts that $H(x,t)$ is a subsolution of \eqref{sish} and
%$w$ is a supersolution of \eqref{equwevol}, as in the proof of
%Theorem \ref{maxpneum} we get
%$$(c(\xs)+\lam)H(\xs,\overline{t})\geq
%(c(\xs)+\lams)w(\xs,\overline{t}),$$ which is a contradiction.

\end{document}